\documentclass[12pt,a4paper]{article}

\usepackage{amsthm}
\usepackage{epsfig}
\usepackage{alltt}
\usepackage{makeidx}
\usepackage{newlfont}
\usepackage{amsmath}
\usepackage{amssymb}
\usepackage{amsfonts}
\usepackage{amscd}
\usepackage[german]{babel}
\usepackage{hyperref}
\usepackage{version}

\makeindex

\newcommand{\SDu}{C'}
\newcommand{\SBu}{\tilde{T}}

\newcommand{\op}{\operatorname}

\begin{document}


\newtheorem{theorem}{Theorem}[section]
\newtheorem{vermutung}[theorem]{Vermutung}
\newtheorem{Theorem}[theorem]{Theorem}
\newtheorem{proposition}[theorem]{Proposition}
\newtheorem{Proposition}[theorem]{Proposition}
\newtheorem{korollar}[theorem]{Korollar}
\newtheorem{Korollar}[theorem]{Korollar}
\newtheorem{lemma}[theorem]{Lemma}
\newtheorem{Lemma}[theorem]{Lemma}
\newcommand{\BLemma}{\begin{lemma}}
\newcommand{\ELemma}{\end{lemma}}
\newtheorem{Satz}[theorem]{Satz}
\newtheorem{problem}[theorem]{Problem}

\theoremstyle{definition}
\newtheorem{definition}[theorem]{Definition}
\newcommand{\BDefinition}{\begin{definition}}
\newcommand{\EDefinition}{\end{definition}}
\newtheorem{Vermutung}[theorem]{Vermutung}
\newtheorem{Definition}[theorem]{Definition}
\theoremstyle{remark}
\newtheorem{Notation}[theorem]{Notation}
\newtheorem{Behauptung}[theorem]{Behauptung}
\newcommand{\BBehauptung}{\begin{Behauptung}}
\newcommand{\EBehauptung}{\end{Behauptung}}
\newtheorem{Bemerkung}[theorem]{Bemerkung} 
\newcommand{\BBemerkung}{\begin{Bemerkung}}
\newcommand{\EBemerkung}{\end{Bemerkung}}
\newcommand{\BProposition}{\begin{Proposition}}
\newcommand{\EProposition}{\end{Proposition}}
\newcommand{\BSatz}{\begin{Satz}}
\newcommand{\ESatz}{\end{Satz}}
\newtheorem{Bemerkungen}[theorem]{Bemerkungen}
\newcommand{\BBemerkungen}{\begin{Bemerkungen}}
\newcommand{\EBemerkungen}{\end{Bemerkungen}}
\newcommand{\BTheorem}{\begin{Theorem}}
\newcommand{\ETheorem}{\end{Theorem}}
\newtheorem{Beispiel}[theorem]{Beispiel} 
\newcommand{\BBeispiel}{\begin{Beispiel}}
\newcommand{\EBeispiel}{\end{Beispiel}}
\newtheorem{Beispiele}[theorem]{Beispiele}
\newcommand{\BBeispiele}{\begin{Beispiele}}
\newcommand{\EBeispiele}{\end{Beispiele}}

\newcommand{\id}{\operatorname{id}}
\newcommand{\EEnd}{\operatorname{End}}
\newcommand{\EExt}{\operatorname{Ext}}
\newcommand{\ext}{\operatorname{ext}}
\newcommand{\Hom}{\operatorname{Hom}}
\newcommand{\Aut}{\operatorname{Aut}}
\newcommand{\RRHom}{\operatorname{RHom}}
\newcommand{\hhom}{\operatorname{hom}}
\newcommand{\HHom}{\operatorname{Hom}}
\newcommand{\res}{\operatorname{res}}
\newcommand{\Res}{\operatorname{Res}}
\newcommand{\Rep}{\operatorname{Rep}}
\newcommand{\resalt}{\operatorname{res}\operatorname{alt}}
\newcommand{\Alt}{\operatorname{Alt}}
\newcommand{\alt}{\operatorname{alt}}
\newcommand{\supp}{\operatorname{supp}}
\newcommand{\card}{\operatorname{card}}
\newcommand{\rk}{\operatorname{rk}}
\newcommand{\coker}{\operatorname{coker}}
\newcommand{\ssum}{\operatorname{sum}}
\newcommand{\SSum}{\operatorname{Sum}}
\newcommand{\um}{{\operatorname{um}}}
\newcommand{\re}{{\operatorname{re}}}
\newcommand{\bsupp}{\operatorname{bsupp}}
\newcommand{\bff}{\operatorname{bf}}
\newcommand{\Rest}{\operatorname{Rest}}
\newcommand{\Quot}{\operatorname{Quot}}
\newcommand{\add}{\operatorname{add}}
\newcommand{\Rg}{\operatorname{Rg}}
\newcommand{\tr}{\operatorname{tr}}

\newcommand{\CA}{{\cal A}}
\newcommand{\CB}{{\cal B}}
\newcommand{\CC}{{\cal C}}
\newcommand{\CE}{{\cal E}}
\newcommand{\CF}{{\cal F}}
\newcommand{\CG}{{\cal G}}
\newcommand{\CH}{{\cal H}}
\newcommand{\CI}{{\cal I}}
\newcommand{\CJ}{{\cal J}}
\newcommand{\CK}{{\cal K}}
\newcommand{\CL}{{\cal L}}
\newcommand{\CM}{{\cal M}}
\newcommand{\CN}{{\cal N}}
\newcommand{\CO}{{\cal O}}
\newcommand{\CP}{{\cal P}}
\newcommand{\CQ}{{\cal Q}}
\newcommand{\CR}{{\cal R}}
\newcommand{\CS}{{\cal S}}
\newcommand{\CT}{{\cal T}}
\newcommand{\CU}{{\cal U}}
\newcommand{\CV}{{\cal V}}
\newcommand{\CW}{{\cal W}}
\newcommand{\CX}{{\cal X}}
\newcommand{\CY}{{\cal Y}}
\newcommand{\CZ}{{\cal Z}}

\newcommand{\ra}{\rightarrow}
\newcommand{\lra}{\longrightarrow}
\newcommand{\thra}{\twoheadrightarrow} 
\newcommand{\sra}{\twoheadrightarrow}
\newcommand{\rsa}{\rightsquigarrow}
\newcommand{\da}{\downarrow}
\newcommand{\dda}{\!\downarrow\!}
\newcommand{\ua}{\uparrow}
\newcommand{\uua}{\!\uparrow\!}
\newcommand{\RA}{\Rightarrow}
\newcommand{\hra}{\hookrightarrow}
\newcommand{\hla}{\hookleftarrow}
\newcommand{\IFF}{\Leftrightarrow}
\newcommand{\U}{U- \mbox {mod}-U}
\newcommand{\sira}{\stackrel{\sim}{\rightarrow}}
\newcommand{\lea}{\leftarrow}

\newcommand{\TA}{{\tilde A}}
\newcommand{\TB}{{\tilde B}}
\newcommand{\TC}{{\tilde C}}
\newcommand{\TD}{{\tilde D}}
\newcommand{\TE}{{\tilde E}}
\newcommand{\TF}{{\tilde F}}
\newcommand{\TG}{{\tilde G}}
\newcommand{\TTH}{{\tilde H}}
\newcommand{\TI}{{\tilde I}}
\newcommand{\TJ}{{\tilde J}}
\newcommand{\TK}{{\tilde K}}
\newcommand{\TL}{{\tilde L}}
\newcommand{\TM}{{\tilde M}}
\newcommand{\TN}{{\tilde N}}
\newcommand{\TO}{{\tilde O}}
\newcommand{\TP}{{\tilde P}}
\newcommand{\TQ}{{\tilde Q}}
\newcommand{\TR}{{\tilde R}}
\newcommand{\TS}{{\tilde S}}
\newcommand{\TT}{{\tilde T}}
\newcommand{\TU}{{\tilde U}}
\newcommand{\TV}{{\tilde V}}
\newcommand{\TW}{{\tilde W}}
\newcommand{\TX}{{\tilde X}}
\newcommand{\TY}{{\tilde Y}}
\newcommand{\TZ}{{\tilde Z}}

\newcommand{\TCA}{{\tilde \CA}}
\newcommand{\TCB}{{\tilde \CB}}
\newcommand{\TCC}{{\tilde \CC}}
\newcommand{\TCE}{{\tilde \CE}}
\newcommand{\TCF}{{\tilde \CF}}
\newcommand{\TCG}{{\tilde \CG}}
\newcommand{\TCH}{{\tilde \CH}}
\newcommand{\TCI}{{\tilde \CI}}
\newcommand{\TCJ}{{\tilde \CJ}}
\newcommand{\TCK}{{\tilde \CK}}
\newcommand{\TCL}{{\tilde \CL}}
\newcommand{\TCM}{{\tilde \CM}}
\newcommand{\TCN}{{\tilde \CN}}
\newcommand{\TCO}{{\tilde \CO}}
\newcommand{\TCP}{{\tilde \CP}}
\newcommand{\TCQ}{{\tilde \CQ}}
\newcommand{\TCR}{{\tilde \CR}}
\newcommand{\TCS}{{\tilde \CS}}
\newcommand{\TCT}{{\tilde \CT}}
\newcommand{\TCU}{{\tilde \CU}}
\newcommand{\TCV}{{\tilde \CV}}
\newcommand{\TCW}{{\tilde \CW}}
\newcommand{\TCX}{{\tilde \CX}}
\newcommand{\TCY}{{\tilde \CY}}
\newcommand{\TCZ}{{\tilde \CZ}}

\newcommand{\td}{\tilde d}

\newcommand{\TDe}{{\tilde \Delta}}
\newcommand{\Tna}{{\tilde \nabla}}
\newcommand{\TTheta}{{\tilde \Theta}}

\newcommand{\reg}{\rangle}
\newcommand{\leg}{\langle}
\newcommand{\la}{\lambda}
\newcommand{\al}{\alpha}
\newcommand{\BMMod}{{\operatorname{-Mod_\DZ-}}}
\newcommand{\mood}{{\operatorname{-mod}}}
\newcommand{\Bmood}{{\operatorname{-mod-}}}
\newcommand{\BMod}{{\operatorname{-Mod-}}}
\newcommand{\MMod}{{\operatorname{-Mod}}}
\newcommand{\Gr}{{\operatorname{Gr}}}
\newcommand{\opp}{{\operatorname{opp}}}
\newcommand{\Ann}{{\operatorname{Ann}}}
\newcommand{\ee}{{\operatorname{ee}}}
\newcommand{\ed}{{\operatorname{ed}}}
\newcommand{\BModbf}{{\operatorname{-Mod^{\bff}_\DZ-}}}

\newcommand{\DA}{{\Bbb A}}
\newcommand{\DC}{{\Bbb C}}
\newcommand{\DP}{{\Bbb P}}
\newcommand{\DDH}{{\Bbb H}}
\newcommand{\DR}{{\Bbb R}}
\newcommand{\DZ}{{\Bbb Z}}
\newcommand{\DN}{{\Bbb N}}
\newcommand{\DQ}{{\Bbb Q}}
\newcommand{\DV}{{\Bbb V}}
\newcommand{\DF}{{\Bbb F}}
\newcommand{\DD}{{\Bbb D}}

\newcommand{\gs}{{\frak g}}
\newcommand{\hs}{{\frak h}}
\newcommand{\bs}{{\frak b}}
\def\ba{\begin{array}}
\def\ea{\end{array}}
\def\bes{\begin{eqnarray*}}
\def\ees{\end{eqnarray*}}
\def\bi{\begin{enumerate}}
\def\ei{\end{enumerate}}
%



\title{Kazhdan-Lusztig-Polynome und unzerlegbare Bimoduln "uber Polynomringen}
\author{Wolfgang Soergel}
\maketitle
\begin{abstract}
Wir entwickeln eine Strategie zum Beweis der Positivit"at
der Koeffizienten von Kazhdan-Lusztig-Polynomen f"ur beliebige
Coxeter-Grup\-pen. Mein Dank gilt Martin H"arterich und
Catharina Stroppel  f"ur 
Korrekturen
zu vorl"aufigen Versionen dieser Arbeit.
Ganz besonders danke ich Parick Polo,
dessen Bemerkungen dazu geholfen haben, die Arbeit noch ganz 
wesentlich zu verbessern.
\end{abstract}

\tableofcontents


\section{Realisierung von Hecke-Algebren durch Bimoduln}

\begin{Notation}
Gegeben eine svelte additive Kategorie $\cal{A}$ bezeichne 
$\langle \cal{A}\rangle$ ihre spaltende
Grothendieck-Gruppe, also die freie abelsche
Gruppe "uber den Objekten mit Relationen $M = M^{\prime} +
M^{\prime\prime}$ wann immer gilt $M \cong M^{\prime} \oplus
M^{\prime\prime}.$
F"ur ein Objekt $M$ bezeichne $\langle M\rangle$ seine Klasse in
$\langle \cal{A}\rangle.$ 
\end{Notation}

\begin{Notation}
Gegeben ein  $\DZ$-graduierter Ring $A$
bezeichne $A\op{-mod}^{\op{f}}_\DZ$ die Kategorie aller endlich
erzeugten
$\DZ$-graduierten $A$-Moduln. 
Wir schreiben $M[n]$ f"ur das Objekt $M$ mit der um $n$
verschobenen $\DZ$-Graduierung, in Formeln $(M[n])_{i} = M_{i+n}$.
\end{Notation}
\begin{Bemerkung}\label{KM}
Gegeben ein K"orper $k$ und 
eine endlich erzeugte kommutative nichtnegativ 
$\DZ$-graduierte ${{k}}$-Algebra $A$ mit $A_0=k$
ist der Endomorphismenring eines endlich erzeugten $\DZ$-graduierten
$A$-Moduls stets endlichdimensional.
Der Endomorphismenring eines unzerlegbaren derartigen Moduls
hat au"ser Null und Eins keine Idempotenten und ist folglich 
unzerlegbar als Rechtsmodul "uber sich selber. Er ist aber auch
der Endomorphismenring dieses unzerlegbaren Rechtsmoduls endlicher
L"ange und folglich lokal.
Mithin gilt in $A\op{-mod}^{\op{f}}_\DZ$ der Satz von Krull-Schmid
mit demselben Beweis wie in \cite{Pierce} 5.4
und die Isomorphieklassen unzerlegbarer Objekte bilden eine Basis
von $\langle A\op{-mod}^{\op{f}}_\DZ\rangle,$
so da"s wir insbesondere haben $\langle N\rangle=\langle M\rangle
\IFF N\cong M.$
\end{Bemerkung}

\begin{Notation}
Sei $(\CW, \CS)$ ein Coxeter-System mit endlich vielen Erzeugern,
$|\CS|<\infty.$ Sei $l :\CW \ra \DN$ die zugeh\"orige
L\"angenfunktion und $\leq$ die Bruhat-Ordnung auf $\CW$. Insbesondere 
bedeutet 
$x<y$ also $x\leq y$, $x\neq y$. 
Auf dem freien $\DZ[v, v^{-1}]$-Modul  
$$\index{$\CH$}\CH=\CH (\CW, \CS) = \bigoplus_{x\in\CW} \DZ[v, v^{-1}] T_{x}$$
\"uber $\CW$ gibt es genau eine Struktur einer assoziativen $\DZ[v, 
v^{-1}]$-Algebra
mit $T_{x} T_{y} = T_{xy}$ falls $l(x)+l(y) = l(xy)$ und $T^{2}_{s} 
=v^{-2}
T_{e} + (v^{-2}-1) T_{s}$ f\"ur alle $s\in \CS,$ siehe \cite{Bou}, IV, \S 
2, 
Exercice 23. Diese assoziative Algebra $\CH$ hei\ss t die Hecke-Algebra 
von 
$(\CW,\CS).$ Sie ist unit"ar mit Eins-Element $T_e,$
wir k"urzen oft $T_e=1$ ab. 
Eigentlich ist es nat"urlicher, mit $q=v^{-2}$ zu arbeiten,
und die ersten Abschnitte dieses Artikels 
w"urden dadurch auch an  Klarheit gewinnen.  
In den sp"ateren Abschnitten
ist es jedoch f"ur eine transparente Darstellung
wichtig, eine Wurzel von $q$ zur Verf"ugung zu haben.
Die Hecke-Algebra 
kann auch beschrieben werden als die assoziative unit"are 
$\DZ[v, v^{-1}]$-Algebra mit
den Erzeugern $\{T_s\}_{s\in\CS},$ den quadratischen 
Relationen $T^{2}_{s} 
=v^{-2}
T_{e} + (v^{-2}-1) T_{s}$ sowie den sogenannten 
Zopfrelationen $T_s T_t \ldots 
T_s=T_tT_s\ldots T_t$ bzw.
$T_s T_t T_{s} \ldots T_t=T_t T_s T_{t}\ldots T_s$ wenn gilt $st\ldots 
s=ts\ldots t$
beziehungsweise $sts\ldots t=tst\ldots s$ f\"ur $s,t\in\CS.$
\end{Notation}

\BDefinition\label{V2}
Sei $(\CW,\CS)$ ein Coxetersystem.
Bezeichne $\CT\subset\CW$ die Menge aller ``Spiegelungen'',
d.h.\ die Menge aller Elemente von $\CW,$ die zu Elementen von $\CS$ 
konjugiert sind. Unter einer {\bf spiegelungstreuen}
Darstellung
unseres Coxetersystems verstehen wir eine Darstellung
${\cal W}\hookrightarrow \op{GL} (V)$ 
in einem endlichdimensionalen Vektorraum "uber einem  
K"orper $k$ der Charakteristik $\op{char} k \neq 2$
mit den folgenden Eigenschaften:
\begin{enumerate}
\item\label{sp1}
Unsere Darstellung ist treu.
\item\label{sp3}
F"ur $x\in\CW$ gilt $\dim (V/V^x)=1\IFF x\in\CT.$ 
In Worten haben genau die Spiegelungen aus $\CW$ 
eine Fixpunktmenge der Kodimension 1 in $V.$ 
\end{enumerate}
\EDefinition
\begin{Bemerkung}
Gegeben eine spiegelungstreue Darstellung  sind die Elemente von $\CT$ 
genau diejenigen
Elemente von $\CW,$ die als Spiegelungen auf $V$ operieren,
die also in anderen Worten $V$ zerlegen in einen eindimensionalen
Eigenraum zum Eigenwert $-1$ und eine Hyperebene aus invarianten
Vektoren. In der Tat kann $x\in \CT$ nicht durch einen unipotenten
Automorphismus operieren, da jeder Automorphismus der Ordnung zwei 
eines Vektorraums "uber einem K"orper der Charakteristik $\neq 2$
diagonalisierbar ist und unsere Bedingung auch ausschlie"st, da"s
$x$ als die Identit"at operiert.
Weiter
kann man die Spiegelungen aus $\CW$ sowohl an ihren
Eigenr"aumen zum Eigenwert $(-1)$ als auch an ihren
Spiegelebenen erkennen,  in Formeln gilt f"ur $t, r \in {\cal T}$ also
$$V^{t}
=  V^{r}\quad \IFF\quad  V^{-t}
=  V^{-r}\quad \IFF\quad t=r.$$
Ist zum Beispiel $G=V^{-t}
=  V^{-r},$ so betrachten wir die kurze exakte Sequenz
$G\hra V\sra V/G.$ Da $t$ und $r$ beide auf
$V/G$ als die Identit"at operieren und
da beim "Ubergang zur transponierten Abbildung die
Dimension der Eigenr"aume gleich bleibt, 
fixiert $tr$ eine Hyperebene. 
Da aber $tr$ auch die
Determinante 1 hat, 
mu"s $tr$ nach Bedingung \ref{sp3} als die Identit"at auf $V$ operieren
und dann nach Bedingung \ref{sp1} die Identit"at in $\CW$ sein.
Wir zeigen im n"achsten Abschnitt, da"s diejenige Darstellung
einer Coxetergruppe, die die Darstellungen der affinen
Weylgruppen auf den Cartan'schen von
Kac-Moody-Algebren verallgemeinert, stets spiegelungstreu ist.
\end{Bemerkung}
\begin{Definition}
Eine Darstellung von $\cal{W},$ bei der Spiegelungen durch Spiegelungen
operieren und bei der man die Spiegelungen an ihren $(-1)$-Eigenr"aumen
unterscheiden kann, nennen wir eine
{\bf spiegelvektortreue} Darstellung. 
\end{Definition}
\begin{Beispiel}
Die geometrische Darstellung einer unendlichen Diedergruppe
ist spiegelvektortreu. Sie ist  aber nicht spiegelungstreu, 
da alle vom neutralen Element verschiedenen Elemente eine 
Fixpunktmenge der Kodimension eins haben, keineswegs nur
die Spiegelungen. 
\end{Beispiel}
\begin{Notation}
Sei nun  zur Vermeidung technischer Notation
$k$ ein unendlicher K"orper und
bezeichne $R = R(V)$ die ${{k}}$-Algebra der regul"aren Funktionen auf
dem einer spiegelvektortreuen Darstellung $V$ zugrundeliegenden
Vektorraum. Wir
versehen $R$ mit einer $\DZ$-Graduierung $R = \bigoplus_{i \in \DZ}
R_{i}$ derart, da"s gilt $R_{2} = V^{\ast}$ und $R_{i} =0$ f"ur
ungerades $i$. Wieder w"are es nat"urlicher, mit der "ublichen Graduierung
zu arbeiten,
und die ersten Abschnitte dieses Artikels 
w"urden dadurch an Klarheit gewinnen.  
Im weiteren Verlauf werden sich unsere Konventionen jedoch bezahlt machen.  
Bezeichne $$\cal{R}=\cal{R}_V\subset R\op{-mod_{\DZ}-}R$$ die Kategorie aller
$\DZ$-graduierten $R$-Bimoduln, die endlich erzeugt sind sowohl als
Rechtsmoduln als auch als Linksmoduln,
und auf denen die $k$-Operationen von rechts und links "ubereinstimmen.
Das Tensorieren $\otimes_R$
"uber $R$ macht $\langle\cal{R}\rangle$ zu
einem Ring.
F"ur $s\in\CS$ bezeichne $R^s\subset R$ den Teilring der $s$-Invarianten.
\end{Notation}

\BTheorem\label{HTT}
Sei $V$ eine spiegelvektortreue Darstellung  eines
Coxetersystems $(\CW,\CS)$ "uber einem
unendlichen K"orper. Bezeichne $\CH$ die Hecke-Algebra und 
$R$ den Ring der polynomialen Funktionen auf $V.$
Es gibt genau einen Ringhomomorphismus
$$\CE :\CH \ra \langle \cal{R}\rangle$$
derart, da"s gilt $\CE (v) = \langle R[1]\rangle$ und $\CE (T_{s}+1) =
\langle R \otimes_{R^{s}} R\rangle \quad \forall s \in \CS.$ 
\ETheorem

\begin{Bemerkung}
Dies Theorem ist eine Variante von
Theorem 1 aus
\cite{HCH}. Mir gef"allt der hier gegebene Beweis besser,
da er keine Kenntnisse
"uber Demazure-Operatoren voraussetzt.
\end{Bemerkung}
\begin{proof}
Die Eindeutigkeit ist klar, da $\CH$ von den $T_{s}+1$ als 
$\DZ[v, v^{-1}]$-Algebra erzeugt wird.
Den Beweis der Existenz m"ussen wir nur im Fall einer Diedergruppe
f"uhren, da ja $\CH$ beschrieben werden kann durch Relationen, 
in denen jeweils nur die Erzeuger zu
zwei einfachen Spiegelungen auftreten.
Wir behandeln diesen Fall in Abschnitt  \ref{DFa}.
\end{proof}
\begin{Definition}
Wir betrachten in der Hecke-Algebra die Elemente $\SBu_x=v^{l(x)}T_x.$
Nach  \cite{KL-C} gibt es genau eine Involution $d:\CH\ra\CH$ 
mit $ d(v)=v^{-1}$ und $d(T_x)=(T_{x^{-1}})^{-1}$ und f"ur jedes $x\in\cal{W}$
genau ein Element $\SDu_{x} \in \cal{H}$ mit $d(\SDu_{x})=\SDu_{x}$ und
$$\SDu_{x} \in \SBu_x+ \sum_y v\DZ[v] \SBu_y$$
Diese Elemente bilden die sogenannte Kazhdan-Lusztig-Basis
der Hecke-Algebra aus \cite{KL-C}.
\end{Definition}
\begin{Vermutung}\label{HV}
Sei $V$ eine spiegelungstreue
Darstellung von $\cal{W}$ "uber einem unendlichen
K"orper und sei $R$ der Ring der regul"aren Funktionen auf $V.$
So gibt es, zumindest im Fall $k = \Bbb{C},$ einen unzerlegbaren
$\Bbb{Z}$-graduierten $R$-Bimodul
$B_{x} \in \cal{R}$ mit $$\CE 
(\SDu_{x}) = \langle B_{x} \rangle.$$
\end{Vermutung}
\begin{Bemerkung}
Ist $k = \Bbb{C}$ und $\CW$ eine endliche Weylgruppe, so wird das
gezeigt in \cite{HCH}.
Ist $\CW$ eine endliche Weylgruppe und $\op{char} k$ mindestens die
Coxeter-Zahl, so wird in \cite{So-R} gezeigt, da"s die Vermutung
"aquivalent ist zu einem Teil der Vermutung von Lusztig zu
Charakteren irreduzibler Darstellungen algebraischer Gruppen "uber
$k$.
F"ur den Fall der Operation der Weylgruppe auf der Cartan'schen einer
affinen Kac-Moody-Algebra wird die
Vermutung bewiesen in \cite{MH}.
\end{Bemerkung}
\begin{Bemerkung}
In \ref{VDI} konstruieren 
wir ein Linksinverses zu $\CE,$ 
dessen explizite Form zeigt, wie aus der vorstehenden Vermutung
die Positivit"at aller Koeffizienten aller Kazhdan-Lusztig-Polynome
folgt. 
\end{Bemerkung}
\begin{Bemerkung}
In  \ref{KL} zeigen wir, da"s
die ``unzerlegbaren Bimoduln im Bild von $\CE$''
zumindest parametrisiert 
werden durch $\CW$. Weiter bestimmen wir in \ref{DHOM} 
die Dimension
der Hom-R"aume zwischen zwei ``Bimoduln im Bild von $\CE$'' und
zeigen, da"s unsere Formeln 
f"ur diese Dimensionen gut mit der Vermutung \ref{HV} zusammenpassen.
\end{Bemerkung}
\begin{Bemerkung}
Noch offen sind insbesondere der Fall universeller Coxetergruppen 
(dort werden die
KL-Polynome von Dyer \cite{Dy} beschrieben), 
und der Fall allgemeiner endlicher
Coxetergruppen, der expliziten Rechnungen mit dem Computer 
zug"anglich sein sollte.
\end{Bemerkung}

\begin{Bemerkung}
Verwandte Untersuchungen findet man 
in \cite{BrM} und (teils unver"offentlichten) Arbeiten
von Dyer \cite{Dyu} und unver"offentlichten Arbeiten
von Peter Fiebig.
\end{Bemerkung}

\section{Eine spiegelungstreue Darstellung}
\begin{Proposition}
Gegeben eine endliche Menge $\cal{S}$ und eine 
Coxetermatrix vom
Typ $\cal{S}$ der Gestalt 
$(m_{s,t})_{s,t \in \cal{S}}$ sei $V$ ein endlichdimensionaler 
reeller Vektorraum und seien gegeben linear 
unabh"angige Vektoren  $(e_{s})_{s\in \cal{S}}$ in $V$
und linear unabh"angige Linearformen
$(e_{s}^{\vee})_{s\in \cal{S}}$ auf $V$ derart, da"s gilt
$$\langle e_{t},e_{s}^{\vee}\rangle 
= -2 \cos (\pi/m_{s,t}) \quad \forall s,t \in \cal{S}.$$
Nehmen wir wie in \cite{Kac} zus"atzlich an, da"s die Dimension
von $V$ kleinstm"oglich ist, 
so liefert die Vorschrift $\rho (s) (v) 
= v - \langle v,e_{s}^{\vee}\rangle e_{s}$ eine spiegelungstreue
Darstellung $\rho:\CW\ra\op{GL}(V)$ unserer Coxetergruppe.  
\end{Proposition}
\begin{proof}
Da"s wir "uberhaupt eine Darstellung erhalten, d.h.\ da"s
die Zopfrelationen erf"ullt sind folgt daraus, da"s
f"ur jedes Paar $s,t$ von einfachen Spiegelungen mit $st$ von
endlicher Ordnung unser Raum zerf"allt in den Schnitt der
Kerne von $e_{s}^{\vee}$ und $e_{t}^{\vee},$ der von $s$ und $t$
"uberhaupt nicht bewegt wird,
und das Erzeugnis von $e_{s}$ und $e_{t},$ auf dem
$s$ und $t$ wie Erzeuger der "ublichen Diedergruppe operieren.

Im allgemeinen ist der von den $e_{s}$ aufgespannte 
Teilraum  $E\subset V$  eine Unterdarstellung, 
die isomorph ist zur ``nat"urlichen Darstellung'' 
aus \cite{Bou}, V.4.3 und die insbesondere
bereits treu ist.
Weiter operiert $\CW$ trivial auf dem simultanen 
Kern $K$ aller $e^{\vee}_{s}$ und 
die nat"urliche Operation auf $(V/K)^{\ast}$ ist wieder isomorph 
zur nat"urlichen Darstellung, d.h.\ die Operation
auf $(V/K)$ ist isomorph zur kontragredienten 
Darstellung $E^\ast$ aus \cite{Bou}, V.4.4.
Unter unserer Minimalit"atsbedingung an die Dimension von $V$ 
haben wir nat"urlich auch $K\subset E$ und damit eine Filtrierung
durch Unterdarstellungen
$$0\subset K\subset E\subset V.$$
Die Paarung $(e_s,e_t)=\frac{1}{2}\langle e_s,e^{\vee}_{t}\rangle$
definiert weiter eine invariante symmetrische  Bilinearform
auf $E$ und eine nichtausgeartete invariante symmetrische  Bilinearform
auf $E/K.$

Sei nun $x \in \cal{W}$ ein Element,
dessen Fixpunktmenge $V^{x} \subset V$ eine Hyperebene ist.
Es gilt zu zeigen $x \in \cal{T}.$
Bezeichne dazu $X \in \op{End} V$ das Bild von $x.$
Wir sind in einem der beiden F"alle $\op{det} X =1$ oder $\op{det} X =-1$
und beginnen mit dem Fall $\op{det} X =1.$ 
W"ahlen wir ein Komplement $\Bbb{R} c$
der Hyperebene $V^{x},$ so hat unser Endomorphismus 
$X$ notwendig die Gestalt
$X (v + \alpha c) = v + \alpha v_{0} + \alpha c$ f"ur 
ein festes $v_0 \in V^{x}.$
Nat"urlich haben wir $V^{x} \supset K.$ Da  unsere Gruppe
treu operiert auf $E$ k"onnen wir $c \in E$ 
w"ahlen und folgern $v_{0} \in E.$
Da auch $E^{\ast} \cong V/K$ eine treue 
Darstellung ist, folgt weiter $v_{0} \not\in
K.$ Jetzt arbeiten wir mit der invarianten 
Bilinearform auf $E$ und folgern
$(v + \alpha c, v + \alpha c) = 
(v + \alpha v_{0} + \alpha c, v + \alpha v_{0} + \alpha c)$ alias
$2 \alpha (v,v_{0}) +
 \alpha^{2} (2(c,v_{0})+(v_{0},v_{0})) = 0$ 
f"ur alle $v \in E$ und $\alpha \in \Bbb{R}.$
Das liefert jedoch $(v,v_{0}) =0$ 
f"ur alle $v\in E$ im Widerspruch zu $v_{0} \not\in K.$

Der Fall $\op{det} X =1$ ist also unm"oglich und 
es bleibt, den Fall $\op{det} X = -1$ zu diskutieren.
Dann operiert $x$ notwendig als Spiegelung 
auf $V$ und auch auf  $V/K \cong E^{\ast}.$
Da $\CW$ treu und transitiv auf den Alkoven
im Tits-Kegel operiert, kann die Fixpunktmenge von $x$ 
keinen Alkoven im Tits-Kegel treffen.
Betrachten wir den fundamentalen dominanten 
Alkoven im Tits-Kegel $C\subset E^{\ast} $ und sein
Bild unter $x$, so werden $C$ und $xC$ nur von 
endlich vielen Spiegelungsebenen $V^z$ zu 
gewissen $z\in \cal{T}$ getrennt.
Andererseits besteht ein geeignetes 
nichtleeres offenes St"uck der Spiegelebene 
von $x$ in $E^\ast$ aus Punkten auf Segmenten, die Punkte von
$C$ mit Punkten von $xC$ verbinden,
man kann hier zum Beispiel das Bild von $C$ unter dem Projektor
$(\op{id}+x)/2$ nehmen.
Dieses St"uck m"u"ste dann von unseren
endlich vielen Spiegelebenen "uberdeckt werden,
und das ist nur m"oglich, wenn die
Spiegelebene von $x$ schon die Spiegelebene einer 
Spiegelung $z$ aus $\cal{T}$ ist.
Dann bilden wir aber $zx$ und finden uns im bereits
ausgeschlossenen Fall $\op{det} =1$
wieder, es sei denn, es gilt $x =z.$
\end{proof}

\section{Weitere Notationen und Formeln}\label{Not}
\begin{Notation}
Gegeben irgendeine endlichdimensionale Darstellung $V$ 
einer Gruppe $\cal{W}$  betrachten wir
f"ur jedes $x \in \CW$  seinen (vertauschten) Graphen
$$ \Gr(x) =\{(x\lambda ,\lambda ) \mid \lambda  \in V\}\subset V\times V$$
und bilden f"ur eine beliebige endliche Teilmenge $A \subset \CW$
in $V\times V$ die Zariski-abgeschlossene Teilmenge
$$ \Gr(A) = \bigcup_{x\in A}  \Gr(x).$$  
\end{Notation}
\begin{Bemerkung}\label{BGr}
F"ur beliebige $y,z\in\CW$ erhalten wir
offensichtlich einen Isomorphismus 
$$ \Gr(y) \cap  \Gr(z) \sira V^{yz^{-1}} $$ 
vermittels der Projektion auf die erste Koordinate und es gilt $$
( \Gr(y) +  \Gr(z))\cap (V\times 0)=\op{im}(yz^{-1}-\op{id})\times 0.$$  
\end{Bemerkung}
\begin{Notation}
Sei nun $k$ ein unendlicher K"orper und
bezeichne wie bisher  $R$ die $k$-Algebra der regul"aren Funktionen auf $V.$
Wir k"urzen stets $\otimes_k=\otimes$ ab.
Identifizieren wir $R\otimes R$ mit den regul"aren Funktionen auf $V\times V$
durch die Vorschrift $(f\otimes g)(\lambda ,\mu)=f(\lambda )g(\mu),$ so bilden 
die regul"aren Funktionen auf $ \Gr(A)$ als Quotient von $R\otimes R$ in
nat"urlicher Weise einen $\DZ$-graduierten $R$-Bimodul. Diesen
$\DZ$-graduierten $R$-Bimodul bezeichnen wir mit
$$R (A) = R ( \Gr(A)) \in R\op{-mod_\DZ-}  R$$
und man pr"uft, da"s er endlich erzeugt ist sowohl als Rechts-
wie auch als Linksmodul "uber $R.$
F"ur $A=\{x_1,\ldots,x_n\}$ schreiben wir auch $ \Gr(A)= \Gr(x_1,\ldots,x_n)$
und
$R(A)=R(x_1,\ldots,x_n)$ oder noch k"urzer $R(x)=R_x,$ $R(x,y)=R_{x,y}.$
Im Fall $A=\{y\mid y\leq x\}$ vereinbaren wir die Abk"urzung 
$R(A)=R(\leq x).$
F"ur die Rechtsoperation von $\CW$ auf $R$ vereinbaren wir eine exponentielle
Schreibweise, also $r^y(\lambda )=r(y\lambda )$ $\forall y\in\CW,$ $\lambda 
\in V,$ $r\in R.$
Bezeichnet $1_y\in R_{y}$ die konstante Funktion $1$ auf $\Gr(y)$ 
alias das Bild von $1\otimes 1$ in $R_y,$ so gilt
dann $r1_y=1_y r^y.$
Die Notationen sind gerade so gew"ahlt, da"s gilt $R(x)\otimes_{R}
R(y) \cong R(xy)$ f"ur alle $x,y \in \CW$.  
\end{Notation}

\section{Der Fall einer Diedergruppe}\label{DFa}

\begin{Bemerkung}\label{DDD}
Ist $\CW$ eine Diedergruppe, so 
haben nach \cite[7.12(a)]{Hu-Re}
die Elemente 
der selbstdualen Basis  die Gestalt
$\SDu_x= v^{l(x)}\sum_{y\leq x} T_{y}.$ 
Wir geben einen Beweis dieser Tatsache als Bemerkung \ref{BD}.  
\end{Bemerkung}

\begin{Theorem}
Sei $(\CW,\CS)$ ein Coxetersystem mit zwei Erzeugern
$|\CS|=2$  und $V$ eine spiegelvektortreue
Darstellung "uber einem unendlichen K"orper.
So ist der Homomorphismus von additiven Gruppen
$$
\begin{array}{cccl}
{\CE} :& \CH& \ra& \langle \cal{R} \rangle\\
&v^{n}\SDu_x &\mapsto&\langle R(\leq x) \rangle[n+l(x)]
\end{array}
$$
ein Ringhomomorphismus.
\end{Theorem}
\begin{Bemerkung}\label{ER}
F"ur jede Spiegelung $s:V\ra V$ induziert
die offensichtliche Abbildung $R\otimes R\sra R(e,s)$ einen
Isomorphismus $R\otimes_{R^s} R\sira R(e,s).$
Um das zu sehen bemerkt man zun"achst, da"s unser
Isomorphismus in spe jedenfalls eine Surjektion ist, und vergleicht dann
die Dimensionen der homogenen Komponenten, wobei man sich st"utzen mag
auf die Eigenraumzerlegung $R=R^s\oplus R^s\alpha$ f"ur $\alpha\in V^\ast$
eine Gleichung der Spiegelebene von $s$ und auf eine kurze
exakte Sequenz $R_s[-2]\hra  R(e,s)\sra R_e.$
Das vorstehende Theorem liefert also in der Tat die Existenzaussage
in Theorem \ref{HTT}.
Der hier wiedergegebene Beweis 
erweitert Argumente  von I. Herrmann \cite{IH}.
\end{Bemerkung}

\begin{proof}
Es reicht offensichtlich, f"ur alle
einfachen Spiegelungen $s$ und alle $x\in\cal{W}$ die Gleichung   
$${\CE}\left((T_s+1)\sum_{y\leq x} T_{y}\right)
={\CE}(T_s+1)\;{\CE}\left(\sum_{y\leq x} T_{y}\right)$$
zu zeigen.
Setzen wir  $A = \{y\in \CW \mid y \leq x\}$, so sind $A
\cup s A$ und $A \cap s A$  von der Gestalt
$\{ y \in \CW \mid y \leq z\}$ f"ur geeignetes $z$ in unserer
Diedergruppe $\CW,$ wenn nicht im zweiten
Fall der Schnitt schlicht leer ist.
Eine explizite Rechnung in der
Hecke-Algebra zeigt nun
$$(T_{s} +1)  \sum_{y \in A} T_{y} = \sum_{y \in A \cup s A}
\!\! T_{y} \;\; +\;\;  v^{-2}\sum_{y \in A \cap s A} \!\! T_{y}.$$
Setzen wir dies Resultat auf der linken Seite ein,
so verwandelt sich unsere Gleichung in die Behauptung eines Isomorphismus
von graduierten Bimoduln
$$
  R (A \cup
s A) \oplus R (A\cap s A)[-2]\cong 
R\otimes_{R^{s}} R\otimes_R  R(A),
$$
der nach einem vorbereitenden Lemma als Proposition \ref{DI} gezeigt wird.
\end{proof}
\begin{Bemerkung}\label{BD}
Um die in \ref{DDD} behaupteten Formeln f"ur die $\SDu_x$
im Diederfall zu zeigen, m"ussen wir nur die Selbstdualit"at der
$\gamma_x=v^{l(x)}\sum_{y\leq x} T_{y}$ nachweisen.
F"ur jede einfache Spiegelung $s$ mit $sx>x$ k"onnen wir jedoch eine
der Gleichungen in unserem Beweis umschreiben zu 
$$v(T_s+1)\gamma_x=
\left\{
  \begin{array}{ll}
\gamma_{sx}+ \gamma_z \text{ mit }z<x&\text{ falls }l(x)>1;\\
\gamma_{sx}&\text{ falls }l(x)\leq 1,
  \end{array}\right.$$
und aus der Selbstdualit"at von $v(T_s+1)$ folgt so induktiv
die Selbstdualit"at aller $\gamma_x.$
\end{Bemerkung}

\BLemma\label{MI}\label{TTT}
Sei $V$ eine endlichdimensionale Darstellung 
einer Gruppe $\cal{W}$ "uber einem 
unendlichen K"orper $k$ einer Charakteristik 
$\neq 2.$
Sei $A \subset \CW$ eine endliche Teilmenge und $s \in \CW$ ein
Element mit $s A = A,$ das auf $V$ als Spiegelung operiert.
So gilt:
\begin{enumerate}
\item
Es gibt einen Isomorphismus von graduierten Bimoduln
$$R\otimes_{R^{s}} R (A) \cong R(A) \oplus R(A) [-2].$$
\item
F"ur $R (A)^{+} \subset R(A)$ die Invarianten unter der Operation
von $s\times \id$  
induziert die Multiplikation einen Isomorphismus
$$ R \otimes_{R^{s}} R(A)^{+} \sira R(A).$$
\end{enumerate}
\ELemma
\begin{proof}
Sei zun"achst $W$ ein beliebiger endlichdimensionaler 
${{k}}$-Vektorraum. Jede Spiegelung $t:W \ra W$ definiert eine Involution
$t : R(W) \ra R(W)$, und w"ahlen wir eine Gleichung $\beta \in
W^{\ast}$ der Spiegelebene $W^{t}$, so k"onnen wir den Operator
$$ \partial_{t} = \partial^{\beta}_{t} : R(W) \ra R(W), \;\;\; f \mapsto
{\frac{f -tf}{2 \beta}}$$
definieren. Ist weiter $X \subset W$ eine Zariski-abgeschlossene
$t$-stabile Teilmenge, so induziert $t$ eine Involution auf $R(X)$
und wir erhalten eine Eigenraumzerlegung der Gestalt $R(X) = R(X)^{+} \oplus
R(X)^{-}$.
Ist keine der irreduziblen Komponenten von $X$ in der Spiegelebene
$W^{t}$ enthalten, so stabilisiert $\partial_t$ den Kern der 
Surjektion $R(W)\sra R(X)$ und induziert folglich
eine Abbildung
$\partial_{t} : R(X) \ra R(X).$ Man sieht, da"s dann 
$\partial_{t}$ und die Multiplikation mit $\beta$ zueinander inverse
Isomorphismen $ R(X)^{+}\cong R(X)^{-}[2]$ von $\DZ$-graduierten
$R^{t}$-Moduln sind.

Sei nun speziell $W = V \times V$ und $t=s \times \id$
f"ur unsere Spiegelung $s\in \CW.$
Wir wenden unsere Erkenntnisse an auf $X =  \Gr(A)$
und erhalten eine Zerlegung $R(A) =
R(A)^{+} \oplus R(A)^{-}$ sowie einen Isomorphismus $ R(A)^{+}
\cong R(A)^{-}[2]$ durch Multiplikation mit $\al\otimes 1=\beta$
f"ur $\al\in V^\ast$ eine Gleichung der Spiegelebene von $s.$
Mit $R = R^{s} \oplus \al R^{s}$ folgt nun Teil 2,
und mit unserer Zerlegung $R(A) =
R(A)^{+} \oplus R(A)^{-}$ und Teil 2 folgt dann auch Teil 1.
\end{proof}

\BProposition\label{DI}
Sei $(\CW,\CS)$ ein Coxetersystem mit zwei Erzeugern
$|\CS|=2$  und $V$ eine spiegelvektortreue
Darstellung von $\CW$ "uber einem unendlichen K"orper.
Seien  $s\in\CS$ und $x\in \CW$
gegeben und bezeichne $A =\{y\in \CW \mid y \leq x\}$
die Menge aller Elemente unterhalb von $x.$
So haben wir in $R\BMMod R$ einen Isomorphismus
$$ R\otimes_{R^{s}} R(A) \cong R(A \cup s A) \oplus R(A\cap s A)
[-2].$$
\EProposition

\begin{proof}[Beweis]
Der Fall $A=sA$ in unserer Proposition hat sich  schon 
durch das vorhergehende Lemma erledigt.
Im Fall $A =\{e\}$ behauptet die Proposition den Isomorphismus
$R(\leq s) \cong R \otimes_{R^{s}} R$ 
aus \ref{ER}.
Damit m"ussen wir nur noch den Fall $A \neq \{e\}$, $A \neq s A$ 
behandeln, also den Fall $A =\{ y\in \CW \mid y \leq x\}$ mit $x \neq e$
und $s x > x$.
Durch explizite Rechnung erhalten wir $A-sA = \{x, rx\}$ f"ur eine von $s$
verschiedene Spiegelung $r \in \CT.$

Nach unseren Annahmen sind die $(-1)$-Eigenr"aume der Spiegelungen aus
$\cal{W}$ paarweise verschieden und spannen eine zweidimensionale
Unterdarstellung $U\subset V$ auf.
Nach  \ref{BGr} gibt es eine  Linearform
$\beta \in V^{\ast} \times V^{\ast},$ 
die auf $\op{Gr} (x) + \op{Gr} (rx)$
verschwindet, nicht aber auf $U\times 0.$ 
Wir betrachten die von 
den Nebenklassen $\overline{\beta}$ und $\overline{1}$ von $\beta$ und $1$
in $R(A)$ erzeugten $R^{s}\otimes R$-Untermoduln $M,N$ und
behaupten
\begin{enumerate}
\item
$M \cong R(A \cap s A)^{+} [-2]$ in $R^{s}\op{-mod_\DZ-} R$;
\item
$N \cong R(A \cup s A)^{+}$ in $R^{s}\op{-mod_\DZ-} R$;
\item
$R(A) = M \oplus N,$
\end{enumerate}
wobei der obere Index $+$ die $(s\times\op{id})$-Invarianten meint.
Sobald das gezeigt ist, folgt die Proposition
aus dem vorhergehenden Lemma \ref{MI}.
Wir m"ussen also nur noch unsere drei Behauptungen zeigen.
\\[2mm]\noindent
1. Um die erste Behauptung einzusehen beachten wir, da"s 
f"ur drei paarweise verschiedene Elemente
$x,y,z \in \cal{W},$ deren L"angen nicht alle dieselbe 
Parit"at haben, stets gilt
$$\op{Gr} (x) + \op{Gr}(y) + \op{Gr}(z)  \supset U\times 0.$$
Um das zu zeigen, d"urfen wir ohne Beschr"ankung der Allgemeinheit $z =e$
und $x,y \in \CT$ annehmen. 
Nach Abschnitt \ref{Not}  gilt dann
$(\op{Gr} (x) + \op{Gr}(e))\cap(U\times 0)=V^{-x}\times 0$
und dasselbe gilt f"ur $y.$  Da wir nun unsere 
Darstellung spiegelvektortreu angenommen hatten, gilt $V^{-x}\neq V^{-y}$
und die Inklusion folgt. 
F"ur $y \not\in \{x,r x\}$ gilt damit $ \Gr(y) +  \Gr(x) +
 \Gr(rx) \supset U \times 0,$ insbesondere verschwindet die Funktion
$\beta$   nicht auf $ \Gr(y)$ f"ur $y \in A \cap sA$.
Ein Element von $R (A)$ annulliert also $\overline{\beta}$ genau
dann, wenn es auf $ \Gr(A\cap s A)$ verschwindet. Wir
folgern, da"s die Multiplikation mit $\overline{\beta}$ einen
Isomorphismus 
$  R(A \cap s A)[-2]\sira R(A)  \overline{\beta}$ liefert.
Das Bild von $R^s\otimes R$ in $R(A \cap s A)$ besteht nun aber genau
aus den $(s\times \op{id})$-Invarianten, folglich schr"ankt unser
Isomorphismus ein zu einem Isomorphismus
$ R(A \cap s A)^+[-2]\sira M.$
\\[2mm]\noindent
2.
Der von $\overline{1}$ in $R(A\cup sA)$ erzeugte $R^s\otimes R$-Untermodul ist 
nat"urlich genau $R(A\cup sA)^+$ und die
Einschr"ankung auf $ \Gr(A)$ liefert eine
Injektion
$$R(A \cup s A)^{+} \hookrightarrow R (A).$$ 
Das zeigt die zweite Behauptung.
\\[2mm]\noindent
3.
Wir zeigen  zun"achst
$R (A) = M +N$.
Bezeichnet $\al\in V^\ast$ eine Gleichung der Spiegelebene $V^s,$ so
gilt nat"urlich $R = R^{s} \oplus \al R^{s}$, also wird $R\otimes
R$ als $R^{s}\otimes R$-Modul erzeugt von $1\otimes 1$ und $\al
\otimes 1$. Es folgt, da"s $R \otimes R$ "uber $R^{s} \otimes R$ auch erzeugt wird von $1\otimes 1$ und beliebiges
nicht unter $(s \times \op{id})$ 
invariantes $\beta \in V^{\ast} \times V^{\ast}.$
Unser $\beta$ kann jedoch nicht 
unter $(s \times \op{id})$ invariant sein, da es sonst auch noch auf
$\op{Gr} (sx)$ verschwinden m"u"ste und damit nach der
Inklusion aus dem Beweis von 1
auf ganz $U\times 0.$
Damit folgt  $R(A) = M+N$.
Wir brauchen nur noch $M \cap N=0$. Das kann man wie
folgt einsehen: F"ur $y \in A \cap s A$ ist die Restriktion
auf $ \Gr(y,sy)=\Gr(y)\cup \Gr(sy)$ 
von jedem Element von $N$ invariant unter
$s \times \id$.
Es reicht zu zeigen, da"s die Restriktion von einem Element von
$M$ auf so ein $\Gr(y,sy)$ nur dann invariant ist unter $s \times \id$,
wenn es auf $ \Gr(y,sy)$ konstant Null ist.
Es reicht dazu, wenn wir zeigen, da"s die Restriktion von
${\beta}$ auf $ \Gr(y,sy)$ nicht invariant ist unter $s
\times \id$ oder auch,
da"s die Restriktion von
${\beta}$ auf $ \Gr(y)+\Gr(sy)$ nicht invariant ist unter $s\times \id.$ 
Die unter $s\times \id$ invarianten Linearformen auf $ \Gr(y)+\Gr(sy)$
m"ussen jedoch verschwinden auf dem Teilraum $V^{-s} \times 0$
und $(\op{Gr}(x) + \op{Gr} (rx))$ schneidet
$U\times 0$  in der davon verschiedenen Gerade $V^{-r}\times 0.$ 
\end{proof}

\section{Von Bimoduln zur"uck zur Hecke-Algebra}
\begin{Notation}
Wir arbeiten von jetzt an stets mit einer
festen spiegelungstreuen Darstellung eines
Coxetersystems "uber einem unendlichen K"orper und
wollen in diesem Abschnitt ein Linksinverses zu unserem
Ringhomomorphismus
$\CE : \CH \ra \langle \cal{R} \rangle$
angeben. Wir beginnen mit einigen Notationen.
Gegeben $B, B^{\prime} \in \cal{R}$ bilden wir $$\HHom
(B, B^{\prime}) = \HHom_{R\otimes R} (B, B^{\prime}) \in
\cal{R}.$$
Hier ist zu verstehen, da"s die Links- bzw. Rechtsoperation von $R$
auf dem $\HHom$-Raum herkommt von der Links- bzw. Rechtsoperation auf $B$ oder
"aquivalent auf $B',$ in Formeln
$(rf)(b)=f(rb)=r(f(b)),$ $(fr)(b)=f(br)=(f(b))r,$ $\forall r\in R,$ $b\in B,$
$f\in\HHom(B,B').$
\end{Notation}
\begin{Notation}
Gegeben ein endlichdimensionaler $\DZ$-graduierter Vektorraum $V =
\bigoplus V_{i}$ definieren wir seine graduierte Dimension
$$\underline{\dim}V = \sum (\dim V_{i}) v^{-i} \in \DZ[v, v^{-1}]$$ und f"ur 
einen
endlich erzeugten $\DZ$-graduierten $R$-Rechtsmodul
$M$ definieren wir seinen graduierten Rang
durch die Vorschrift
$$\underline{\rk} M = \underline{\dim} (M/ MR_{>0})\in \DZ[v, v^{-1}].$$
Es gilt also $\underline{\dim} (V[1]) = v (\underline{\dim} V)$ und 
$\underline{\rk} (M[1])= v (\underline{\rk} M)$.
Nat"urlich ist unser graduierter Rang nur f"ur freie Moduln ein
vern"unftiger Begriff, und wir werden ihn auch nur f"ur freie
Moduln verwenden. Mit $\overline{\underline{\rk}} M$ bezeichnen wir 
das Bild von $\underline{\rk} M$ unter der Substitution
$v \mapsto v^{-1}.$ 
Das folgende Theorem motiviert den ganzen 
Abschnitt.   
\end{Notation}
\begin{Theorem}\label{VBI}\label{VDI}
Unsere Abbildung
$\CE : \CH \ra \langle \cal{R}\rangle$ hat als
Linksinverse die Abbildung $\langle \cal{R}\rangle
\ra \CH$ gegeben durch die Vorschrift
$$\langle B \rangle \mapsto \sum_{x\in\CW} \overline{\underline{\rk}}
\HHom (B, R_x) \; T_{x}.$$
\end{Theorem}

\begin{proof}
Dies Theorem ist eine direkte Konsequenz
aus \ref{HD} und \ref{KOO}, die im
Folgenden ohne R"uckgriff auf das Theorem bewiesen werden.
\end{proof}
\begin{Definition}
Wir 
definieren  f"ur jeden $R$-Bimodul $B$ und jede
Teilmenge $A \subset \CW$ den Unterbimodul
$$\Gamma_{A} B =\{ b \in B \mid \supp b \subset  \Gr(A)\}$$
aller Elemente mit Tr"ager  in $ \Gr(A) .$ 
Wir setzen $\Gamma_{\geq i} B=\Gamma_{\{x\in\CW\mid l(x)\geq i\}} B$
und definieren die Kategorie $$\CF_\Delta\subset \cal{R}$$ als die volle 
Unterkategorie aller graduierten Bimoduln $B\in \cal{R}$ derart, da"s 
    $B$ Tr"ager in einer Menge der Gestalt $\op{Gr}(A)$ hat
f"ur endliches $A\subset\CW$ und
   da"s der
Quotient $\Gamma_{\geq i} B/\Gamma_{\geq i+1}B$ 
f"ur alle $i$ isomorph ist zu einer endlichen direkten Summe
von graduierten Bimoduln der Form $R_{x}[\nu]$ mit $l(x)=i$ und $\nu\in\DZ.$
\end{Definition}
\begin{Bemerkung}
 Sicher ist $B\mapsto 
\Gamma_{\geq i} B/\Gamma_{\geq i+1}B$ additiv in $B,$ 
folglich ist $\CF_\Delta$ stabil unter dem 
Bilden von endlichen direkten Summen
und mit Krull-Schmid \ref{KM}  auch unter dem 
Bilden von direkten Summanden.
\end{Bemerkung}
\begin{Notation}
Es ist nun bequem, die graduierten Bimoduln $$\Delta_{x}=R_{x}[-l(x)]$$ 
einzuf"uhren und in der Hecke-Algebra mit $\SBu_x=v^{l(x)}T_x$ zu arbeiten. 
F"ur die 
Vielfachheit des Summanden $\Delta_{x}[\nu]$ in einer 
und jeder Zerlegung 
von $\Gamma_{\geq i} B/\Gamma_{\geq i+1}B$ f"ur $i=l(x)$ 
vereinbaren
wir die Notation $(B:\Delta_{x}[\nu]).$ 
Schlie"slich vereinbaren wir die 
Abk"urzung $R[1]\otimes_{R^{s}} M=\theta_s M.$  
\end{Notation}

\begin{Proposition}\label{ZD}\label{HD}
Sei $s \in \CS$ eine einfache Spiegelung.
\begin{enumerate}
\item
Mit $B$ liegt auch $R\otimes_{R^{s}} B$ in $\CF_{\Delta}$.
\item
Definieren wir die Abbildung
$h_{\Delta}: \CF_{\Delta} \ra \CH$ durch die Vorschrift $B
\mapsto \sum_{x,\nu}   (B:\Delta_{x}[\nu]) v^{\nu} \SBu_{x}$, so
kommutieren f"ur alle $s\in\CS$ die beiden Diagramme
$$\begin{CD}
\CF_{\Delta} @>>> \CH\\
@V\theta_sVV    @VV(\SBu_{s}+v)\cdot V\\
\CF_{\Delta}@>>> \CH
\end{CD}\hspace{2.5cm}
\begin{CD}
\CF_{\Delta} @>>> \CH\\
@V[1]VV    @VVv\cdot V\\
\CF_{\Delta}@>>> \CH
\end{CD}$$
\item\label{DR}
Unsere Abbildung $\CE$ faktorisiert "uber eine Abbildung
$\CE:\CH\ra \langle\CF_{\Delta}\rangle$ und unser 
$h_\Delta: \langle\CF_{\Delta}\rangle\ra\CH$ ist zu diesem $\CE$
linksinvers.
\end{enumerate}
\end{Proposition}
\begin{proof}
Zum Beweis dieser Proposition ben"otigen wir
\BLemma\label{Ext}
Sei $W$ ein endlichdimensionaler Vektorraum und seien $U,V \subset
W$ zwei affine Teilr"aume.
Nur dann ist $\EExt^{1}_{R(W)} (R(U), R(V))$ verschieden von Null,
wenn gilt $V \cap U=V$ oder wenn $V \cap U$ eine
Hyperebene in $V$ ist. 
Im letzteren Fall ist $\EExt^{1}_{R(W)}
(R(U),R(V))$ ein freier $R(U\cap V)$-Modul vom Rang Eins, erzeugt von
der Klasse einer beliebigen kurzen exakten Sequenz der Form
$$R (V)[-2] \overset{\al}{\hookrightarrow} R(U\cup V)
\twoheadrightarrow R(U)$$
f"ur $\al \in W^{\ast}$ mit $\al|_U = 0$, $\al|_V
\neq  0.$
\ELemma
\begin{proof}
Sind $F$ und $G$ freie Moduln von endlichem Rang "uber den
${{k}}$-Algebren $A$ und $B$ und sind $M$ und $N$ irgendwelche
Moduln "uber $A$ und $B$, so gilt 
offensichtlich (siehe \cite{BouA1}, I \S 4 und II \S 11)
$$\HHom_{A} (F, M) \otimes \HHom_{B} (G,N) \sira \HHom_{A \otimes
B}
(F \otimes G, M \otimes N).$$
Sind unsere Algebren noethersch und $M^{\prime}$ bzw.\
$N^{\prime}$ endlich erzeugte Moduln "uber $A$ bzw.\ $B$, so
finden wir Aufl"osungen $F^{\bullet} \twoheadrightarrow
M^{\prime}$ bzw.\ $G^{\bullet} \twoheadrightarrow N^{\prime}$ mit
$F^{i}$ bzw.\ $G^{j}$ frei von endlichem Rang "uber $A$ bzw.\
$B$.
Dann wird $F^{\bullet} \otimes G^{\bullet}$ eine freie Aufl"osung
von $M^{\prime} \otimes N^{\prime}$ und wir folgern
$$\begin{array}{ccl}
\EExt^{n}_{A\otimes B} (M^{\prime} \otimes N^{\prime}, M\otimes
N)&=& H^{n}\HHom_{A\otimes B} (F^{\bullet} \otimes G^{\bullet}, M
\otimes N)\\
&=&H^{n} (\HHom_{A} (F^{\bullet},M) \otimes
\HHom_{B}(G^{\bullet},N))\\
&=& \bigoplus_{i+j=n} \EExt^{i}_{A} (M^{\prime}, M) \otimes
\EExt^{j}_{B} (N^{\prime},N)
\end{array}$$
Wir k"onnen uns so im Lemma auf die F"alle mit $\dim_{{{k}}} W =1$
zur"uckziehen, und diese behandelt man leicht explizit.
Genauer d"urfen wir 
ohne Beschr"ankung der Allgemeinheit $U$ und $V$ linear
annehmen und finden eine Zerlegung 
$W = S \oplus U^{\prime} \oplus V^{\prime} \oplus W^{\prime}$ mit $U = S \oplus
U^{\prime}$ und $V = S \oplus V^{\prime}.$
Notieren wir die Dimensionen $\op{dim}W =  s + u +v +w,$ so erhalten wir
$\op{Ext}^{\bullet}_{R(W)} (R (U), R(V)) 
\cong \op{Ext}^{\bullet}_{k[X]} (k[X], k[X])^{s}
\otimes \op{Ext}^{\bullet}_{k[X]} (k[X],k)^{u}
\otimes \op{Ext}^{\bullet}_{k[X]}(k,k[X])^{v} 
\otimes \op{Ext}^{\bullet}_{k[X]} (k,k)^{w}.$
Damit kann es nur im Fall $v \leq 1$ "uberhaupt $\op{Ext}^{1}$ 
geben und im
Fall $v =1$ ist $\op{Ext}^{1}$ gerade $k[X]^{s}.$
\end{proof}\noindent
Jetzt zeigen wir die Proposition.
Sei $s \in \CS$ unsere feste einfache Spiegelung. 
Sicher k"onnen wir unsere Filtrierung $\Gamma_{\geq i}$ von $B$ 
verfeinern durch
gewisse $\Gamma_{\geq j}B$ mit $ j \in \DZ+0,5$ derart, da"s
f"ur alle $i\in\DZ$ die Quotienten 
$\Gamma_{\geq i}B/\Gamma_{\geq i+0,5}B$ bzw.\
$\Gamma_{\geq i-0,5}B/\Gamma_{\geq i}B$ Summen sind von $R (x) [\nu]$ mit
$x > sx$ bzw.\ $ x < sx$.
Mit diesen Wahlen unterscheiden sich die Parameter $x,y$ zweier
m"oglicher Subquotienten von $\Gamma_{\geq i-0,5}B/\Gamma_{\geq i+0,5}B$ nur
dann um eine Spiegelung, wenn gilt $y=sx$:
In der Tat haben wir f"ur eine beliebige Spiegelung $t$ stets
$t x > x $ oder $t x < x$, und aus $sy >y < x>s x$ folgt $y
\leq sx$ nach einem Resultat von Deodhar, der sogenannten
Eigenschaft Z von Coxetergruppen.

Alle Erweiterungen in $\EExt^{1}_{R\otimes R} (R_{x}, R_{sx})$
sterben aber unter der Restriktion auf $R^{s}\otimes R.$
In der Tat spaltet die Restriktion $R_{x,sx}\sra R_{x}$ "uber
$R^s\otimes R,$ da sie in den Notationen von \ref{TTT} einen
Isomorphismus $R_{x,sx}^+\sira R_{x}$ induziert.
Nach Lemma \ref{Ext}
ist damit die Einschr"ankung auf $R^{s}\otimes R$ von
$\Gamma_{\geq i-0,5}B/\Gamma_{\geq i+0,5}B$ isomorph zu einer direkten Summe
von Kopien von gewissen $R_x [\nu]$ mit $x>sx$ und $l(x)=i,$ 
und diese treten auf mit Vielfachheit
$(B:R_{x}[\nu]) + (B:R_{sx}[\nu]).$
Tensorieren wir nun wieder hoch vermittels $R\otimes_{R^{s}}$ und
beachten die kurzen exakten Sequenzen $R_{x} [-2] \hookrightarrow
R\otimes_{R^{s}} R_{x} \twoheadrightarrow R_{sx}$, so erhalten wir
die erste Behauptung sowie (immer unter der Annahme $x>sx)$ die Formeln
$$\begin{array}{llll}
(R\otimes_{R^{s}} B:R_{x}[\nu]) &=&
(B:R_{x}[\nu+2])&+\; (B:R_{sx}[\nu+2])\\
(R\otimes_{R^{s}} B:R_{sx}[\nu]) &=&
(B:R_{x}[\nu])&+\;(B:R_{sx}[\nu])\end{array}$$
Das k"onnen wir umschreiben zu
$$\begin{array}{llll}
(\theta_s B : \Delta_{x} [\nu])& = &(B: \Delta_{x} [\nu
+ 1]) &+\; (B: \Delta_{sx} [\nu])\\
(R[1] \otimes_{R^{s}} B: \Delta_{sx} [\nu]) & = & (B : \Delta_{x}
[\nu]) &+\; (B : \Delta_{sx} [\nu -1])
\end{array}$$
Schreiben wir nun ein Element $H\in\CH$ nur f"ur diesen Beweis 
in der Form $H=\sum_{x,\nu}(H:v^{\nu}\SBu_x) \;v^{\nu}\SBu_x,$
so gilt in der Hecke-Algebra ganz analog
$$\begin{array}{llll}
((\SBu_s+v)H:v^{\nu}\SBu_{x}) &=&
(H:v^{\nu+1}\SBu_{x})&+\;(H:v^{\nu}\SBu_{sx})\\
((\SBu_s+v)H:v^{\nu}\SBu_{sx}) &=&
(H:v^{\nu}\SBu_{x})&+\;(H:v^{\nu-1}\SBu_{sx})
\end{array}$$
Das zeigt die zweite Behauptung.
Um die letzte Behauptung zu zeigen, bemerken wir zun"achst,
da"s wir $\langle\cal{R}\rangle$ vermittels 
des Ringhomomorphismus $\cal{E}$ aus \ref{HTT} insbesondere
als $\cal{H}$-Linksmodul auffassen k"onnen.
Der erste Teil der Proposition sagt dann, da"s 
$\langle\cal{F}_\Delta\rangle\subset \langle\cal{R}\rangle$
ein $\cal{H}$-Untermodul ist und der Zweite,
da"s die Abbildung $h_\Delta:\langle\cal{F}_\Delta\rangle
\ra \CH$ ein Homomorphismus von $\CH$-Moduln ist.
Da $\cal{E}(1)=\langle R_e\rangle$ zu $\cal{F}_\Delta$ geh"ort,
faktorisiert $\cal{E}$ damit in der Tat "uber $\langle\cal{F}_\Delta\rangle,$
und da $h_\Delta\circ\CE:\CH\ra\CH$ ein Homomorphismus von
$\CH$-Moduln ist, der die Eins auf die Eins wirft, mu"s diese
Verkn"upfung bereits die Identit"at sein.
\end{proof}
Um Theorem \ref{VDI} zu zeigen, brauchen wir auch einen ``dualen'' Zugang.
Genauer betrachten wir 
die Filtrierung unserer Bimoduln durch die
$\Gamma_{\leq i} B=\Gamma_{\{x\in\CW\mid l(x)\leq i\}} B$
und definieren die Kategorie $$\CF_\nabla\subset \cal{R}$$ als die volle 
Unterkategorie aller graduierten Bimoduln $B\in \cal{R}$ derart, da"s 
ihr Tr"ager in $\op{Gr}(A)$ liegt f"ur endliches $A\subset \CW$ und da"s
die Quotienten f"ur alle $i$ isomorph sind zu einer endlichen direkten Summe
von graduierten Bimoduln der Form $R_{x}[\nu]$ mit $l(x)=i$ und $\nu\in\DZ.$
In diesem Zusammenhang ist es bequem und
nat"urlich, mit $$\nabla_{x}=R_{x}[l(x)]$$ zu arbeiten. 
F"ur die Vielfachheit von $\nabla_{x}[\nu]$ in einer Zerlegung 
in eine direkte Summe von
$\Gamma_{\leq l(x)} B/\Gamma_{\leq l(x)-1}B$ 
vereinbaren wir die Notation $(B:\nabla_{x}[\nu])$.

Man beachte, da"s sich die Multiplizit"aten 
$(B:\nabla_{x}[\nu])$ und $(B:\Delta_{x}[\nu])$
auf Subquotienten in verschiedenen Filtrierungen
beziehen. Auch f"ur 
$B\in\CF_\Delta\cap \CF_\nabla$ sind also die Multiplizit"aten 
$(B:\Delta_{x}[l(x)+\nu])$
und $(B:\nabla_{x}[-l(x)+\nu])$ im allgemeinen verschieden.
Ganz analog wie eben haben wir nun
\begin{Proposition}\label{ZN}\label{Nabla}
Sei $s \in \CS$ eine einfache Spiegelung.
\begin{enumerate}
\item
Mit $B$ liegt auch $R\otimes_{R^{s}} B$ in $\CF_{\nabla}$.
\item
Definieren wir die Abbildung
$h_{\nabla}: \CF_{\nabla} \ra \CH$ durch die Vorschrift $B
\mapsto \sum_{x,\mu}   (B:\nabla_{x}[\mu])\; v^{-\mu} \SBu_{x}$, so
kommutieren f"ur alle $s\in \CS$ die beiden Diagramme
$$\begin{CD}
\CF_{\nabla} @>>> \CH\\
@V\theta_sVV    @VV(\SBu_{s}+v)\cdot V\\
\CF_{\nabla}@>>> \CH
\end{CD}\hspace{2.5cm}
\begin{CD}
\CF_{\nabla} @>>> \CH\\
@V[1]VV    @VVv^{-1}\cdot V\\
\CF_{\nabla}@>>> \CH
\end{CD}$$
\item
Die Verkn"upfung $d\circ h_\nabla$ ist linksinvers 
zu $\CE:\CH\ra \langle\CF_\nabla\rangle.$
\end{enumerate}
\end{Proposition}
\begin{proof}
Wir betrachten den Funktor $$D = \HHom_{-R}(\;,
R):  \cal{R}\ra R\op{-mod_\DZ-}R,$$
wo wir unseren Raum von Homomorphismen von $R$-Rechtsmoduln
versehen mit der offensichtlichen $\DZ$-Graduierung
und die Rechts- bzw. Linksoperation auf $DB$ definieren vermittels
der Rechts- bzw. Linksoperation auf einem Bimodul $B,$ in Formeln
$(rf)(b)=f(rb)$ und $(fr)(b)=f(br)$ f"ur alle $b\in B,$ $r\in R,$ $f\in DB.$
Ich will es vermeiden, im allgemeinen zu diskutieren, ob ein
Objekt aus $\CR$ unter $D$ wieder in $\CR$ landet. In jedem Fall
haben wir $DR_x\cong R_x$ und $D(M[\nu])=(DM)[-\nu],$ 
unser Funktor induziert  eine "Aquivalenz von Kategorien $D :
\CF_{\nabla} \sira \CF_{\Delta}^{\opp}$ und es gilt offensichtlich
$h_{\nabla} = h_{\Delta} \circ D$.
Es bleibt also, f"ur alle $M \in \CF_{\nabla}$ einen Isomorphismus
$\theta_s DM \cong D \theta_s M$
nachzuweisen. Dazu zeigen wir erst einmal
\begin{Proposition}\label{SA}
\begin{enumerate}
\item
Die Funktoren von $R^{s}\op{-mod}_{\DZ}$ nach $R\op{-mod}_{\DZ}$
mit $M \mapsto R[2] \otimes_{R^{s}} M$ und $M \mapsto
\HHom_{R^{s}} (R,M)$ 
sind nat"urlich "aquivalent.
\item
Der Funktor $R[1]\otimes_{R^s}:R\op{-mod}_{\DZ}\ra R\op{-mod}_{\DZ}$ ist 
selbstadjungiert.
\end{enumerate}
\end{Proposition}
\begin{proof}
Da $R$ frei ist von endlichem Rang als $R^{s}$-Modul, l"a"st
sich unser Hom-Funktor auch schreiben in der Form
$$\HHom_{R^{s}} (R,-) = \HHom_{R^{s}}
(R,R^{s})\otimes_{R^{s}}-.$$
Nun ist genauer $R$ frei "uber $R^{s}$ mit Basis $1,\al$ f"ur
$\al$ eine Gleichung der Spiegelebene, $s(\al) = -\al$.
Die duale Basis von $\HHom_{R^{s}} (R,R^{s}) $ "uber $R^{s}$
notieren wir $1^{\ast}, \al^{\ast}.$ Die Multiplikation mit
$\al \in R$ wird in dieser Basis 
gegeben durch $\al \al^{\ast} = 1^{\ast}$, $\al
1^{\ast} = \al^{2} \al^{\ast}$.
Die Wahl von $\al$ definiert folglich einen Isomorphismus von
$R$-Moduln
$$R[2]\sira \HHom_{R^{s}} (R, R^{s}), 1 \mapsto \al^{\ast}, \al
\mapsto 1^{\ast},$$
und die erste Behauptung ist bewiesen.
Die zweite Behauptung folgt, da ja unsere beiden Funktoren aus
der ersten Behauptung 
bis auf eine Verschiebung der Graduierung
gerade die beiden Adjungierten der Einschr"ankung auf
$R^s$ sind.
\end{proof}\noindent
Damit erhalten wir in der Tat Isomorphismen
$$\begin{array}{ccl}
\theta_s (DM) & \cong&R[1]\otimes_{R^{s}} DM\\
&\cong & \HHom_{R^{s}}
(R[1],\HHom_{-R} (M,R))\\
&\cong &\HHom_{- R} (\theta_s M, R)\\
& \cong &D (\theta_s M)\end{array}$$
Der dritte Teil der Proposition \ref{Nabla} folgt wieder daraus,
da"s wegen $d(\SBu_s+v)=\SBu_s+v$ die Verkn"upfung
$d\circ h_\nabla\circ\CE:\CH\ra\CH$ ein Homomorphismus von 
$\CH$-Linksmoduln
ist mit $1\mapsto 1.$
\end{proof}
\begin{definition}
Wir definieren nun die Kategorie $\CB\subset \cal{R}$
als die Kategorie aller graduierten
Bimoduln $B\in\cal{R},$ deren Klasse $\langle B\rangle$
im Bild unseres Morphismus  
$\CE:\CH\ra \langle \cal{R}\rangle$ liegt, und nennen
die Objekte von $\CB$ die {\bf speziellen Bimoduln}. 
\end{definition}
\begin{Bemerkung}
Sicher ist $\CB$ stabil unter endlichen direkten Summen und unter 
Verschiebungen
der Graduierung. Erst sp"ater in \ref{DSS} werden wir zeigen
k"onnen, da"s $\CB$ auch stabil ist unter dem Bilden von
direkten Summanden.
Um ein Kriterium daf"ur zu erhalten, wann ein
Bimodul $B$ zu $\cal{B}$ geh"ort, betrachten wir
zun"achst f"ur eine beliebige endliche Sequenz
$\underline{s}=(r,\ldots,t)$ von einfachen Spiegelungen
in
der Hecke-Algebra
das Element $b(\underline{s})=(T_r+1)\ldots (T_t+1)$  und bilden den Bimodul
$$B(\underline{s})=R\otimes_{R^r}\ldots R\otimes _{R^t}R.$$
\end{Bemerkung}
\begin{Lemma}\label{LN}
Ein graduierter Bimodul $B\in \cal{R}$ ist speziell genau dann, wenn 
es Objekte $C,D\in \cal{R}$ gibt, die jeweils endliche direkte  Summen
sind von Objekten der Gestalt $B(\underline{s})[n]$ und so, da"s gilt
$$B\oplus C\cong D$$
\end{Lemma}
\begin{proof}
Da gilt
$\langle B(\underline{s})[n]\rangle=\cal{E}(v^nb(\underline{s}))$ 
sind die $B(\underline{s})[n]$ speziell.
Das zeigt, da"s unser Kriterium hinreichend ist.
Da die $v^nb(\underline{s})$ schon $\CH$ erzeugen als 
abelsche Gruppe, ist es auch notwendig. 
\end{proof}
\begin{Bemerkung}
Insbesondere faktorisiert $\CE:\CH\ra \langle \cal{R}\rangle$ "uber 
die spaltende Grothendieckgruppe 
$\langle\CB \rangle$  
der additiven Kategorie $\CB$
und aus \ref{ZD} und \ref{ZN} folgt 
$\CB\subset\CF_\Delta\cap\CF_\nabla.$ 
\end{Bemerkung}

\begin{Theorem}\label{HOM}\label{DHOM}
F"ur $M \in \CF_{\Delta}$, $ N \in \CB$ 
und desgleichen f"ur $M\in\CB, $ $N\in \CF_\nabla$ ist $\HHom
(M,N)$ graduiert frei als $R$-Rechtsmodul vom Rang 
$${\underline{\rk}}\HHom
(M,N)=
\sum_{x,\nu,\mu} (M: \Delta_{x} [\nu])(N:\nabla_{x} [\mu])v^{\mu
-\nu}.$$
\end{Theorem}
\begin{proof}
Wir behandeln nur den Fall $M \in \CF_{\Delta}$, $N \in \CB$,
der andere Fall geht genauso.
Bezeichne $i :\CH \ra \CH$ die Antiinvolution mit $i(v) = v$, $i
(\SBu_{x})=\SBu_{x^{-1}}$.
Wir betrachten die symmetrische $\DZ[v, v^{-1}]$-bilineare Paarung $\langle
\; , \; \rangle : \CH \times \CH \ra \DZ[v, v^{-1}]$, die jedem Paar $(F,G)$ 
den
Koeffizienten von $\SBu_{e} = T_{e}$ in der Darstellung von
$i(F) G$ 
als Linearkombination der $T$ oder "aquivalent der $\SBu$ zuordnet.
Sie kann auch explizit beschrieben werden durch $\langle \SBu_{x},
\SBu_{y} \rangle = \delta_{xy},$ vergleiche \cite{Lu-AfC},\S1.4.
Die gesuchte Formel erh"alt mit
dieser Notation die Gestalt
$$\overline{\underline{\rk}} \HHom (M,N) = \langle h_{\Delta} M,
h_{\nabla} N \rangle.$$
Nach \ref{LN} d"urfen wir uns hier zun"achst einmal auf den
Fall $N=B(\underline{s})$ zur"uckziehen.
Wir erkennen weiter mithilfe von \ref{SA}, da"s unsere Formel 
richtig ist f"ur das Paar $(R\otimes_{R^s}M ,N)$
genau dann, wenn sie richtig ist f"ur das Paar $(M,R\otimes_{R^s}N ).$
Ohne alle Probleme sehen wir auch, 
da"s unsere Formel 
richtig ist f"ur das Paar $(M[1],N)$
genau dann, wenn sie richtig ist f"ur das Paar $(M,N[-1]).$
Mit diesen Erkenntnissen k"onnen wir uns dann 
schlie"slich sogar auf den
Fall $N=R_e$ zur"uckziehen, und der ist offensichtlich.
\end{proof}
\begin{korollar}\label{KOO}
Sei $B\in\CB$ einer unserer speziellen Bimoduln.
So gilt
$$\begin{array}{lll}
h_\Delta(B)&=& \sum_{x\in\CW} \overline{\underline{\rk}}
\HHom (B, R_{x}) \; T_{x},\\[2mm]
h_\nabla(B)&=& \sum_{x\in\CW} \overline{\underline{\rk}}
\HHom (R_{x},B) \; T_{x}.
\end{array}$$
\end{korollar}
\begin{proof}
Nach dem vorhergehenden Theorem \ref{HOM} haben wir
$$\begin{array}{lll}
\underline{\rk}\HHom (B, R_{x})&=&\underline{\rk}\HHom (B, \nabla_{x}[-l(x)])\\
&=&\sum_{\nu}(B:\Delta_{x}[\nu])v^{-l(x)-\nu}\\[3mm]
\underline{\rk}\HHom (R_{x},B)&=&\underline{\rk}\HHom (\Delta_{x}[l(x)],B)\\
&=&\sum_{\mu}(B:\nabla_{x}[\mu])v^{-l(x)+\mu}
\end{array}$$ 
und das Korollar folgt damit aus den Definitionen.
\end{proof}

\section{Klassifikation der unzerlegbaren speziellen Bimoduln}
\begin{Notation}
Wir arbeiten weiter mit einer spiegelungstreuen Darstellung
"uber einem unendlichen K"orper. F"ur $B\in\CB$ und $y\in \CW$ k"urzen wir ab
$\Gamma_{\leq y}B/\Gamma_{< y} B=\Gamma^\leq_y  B,$ ebenso auch 
$\Gamma_{\geq y}B/\Gamma_{> y} B=\Gamma^\geq_y  B$ und 
$B/\Gamma_{\neq y}B=\Gamma^y B.$ 
\end{Notation}

\begin{Bemerkung}\label{BB}
Gegeben ein Bimodul mit einer 
endlichen Filtrierung 
durch Unterbimoduln mit zu geeigneten $R_x$ isomorphen
Subquotienten ist der Tr"ager
jedes Schnittes $s \in B$ eine Vereinigung von Graphen von Elementen
unserer Coxetergruppe. Um das zu zeigen, d"urfen wir ohne Beschr"ankung der
Allgemeinheit $\op{supp} s$ irreduzibel annehmen.
Nun liefert aber jedes von Null verschiedene 
$s \in B$ ein von Null verschiedenes
Element $\bar{s}$ in einem geeigneten 
Subquotienten der $\Delta$-Fahne, und
es folgt erst $\op{supp} s \supset \op{supp} \bar{s} = \op{Gr} (x)$ f"ur
geeignetes $x \in W$ und dann wegen der Irreduzibilit"at von $\op{supp}s$
sogar $\op{supp}s = \op{Gr}(x).$
 Insbesondere haben wir 
$\Gamma_{y}B\cap \Gamma_{\neq y}B=0$
und damit Einbettungen $\Gamma_{y}B\hra \Gamma^{\leq}_{y} B$ sowie 
$\Gamma_y^{\geq}B\hra \Gamma^{y}B.$  
Unter denselben Voraussetzungen an $B$ "andert auch die Multiplikation
von rechts mit einem von Null verschiedenen Element $p\in R$ den
Tr"ager nicht. Wir haben also
$(\Gamma_A B)p= \Gamma_A (Bp)$ und dann auch
$(\Gamma_y^{\geq}B)p\cong\Gamma_y^{\geq}(Bp)$ und dergleichen.
Wir werden uns diese Identifikationen zunutze machen,
um die Klammern bei derartigen Ausdr"ucken wegzulassen.
\end{Bemerkung}

\BLemma\label{PP}
Sei 
$B\in\cal{F}_\nabla$ gegeben und
sei $y_{0}, y_{1},y_{2}, \ldots $ eine Aufz"ahlung der
Elemente von $\CW$ derart, da"s in der Bruhat-Ordnung gr"o"sere
Elemente stets auch einen gr"o"seren Index haben. Bezeichne
$C(k) = \{y_{0}, \ldots, y_{k}\}$ die Menge der 
ersten $k$ Elemente in unserer Aufz"ahlung und sei $y_k=y$.
So induziert die offensichtliche Abbildung einen Isomorphismus
$$\Gamma_{\leq y} B/ \Gamma_{<y} B\sira \Gamma_{C(k)} B/\Gamma_{C(k-1)} B,$$
beide Seiten sind direkte Summen von Objekten der Gestalt $\nabla_{y}
[\nu],$
und ein
$\nabla_{y}
[\nu]$ kommt in diesem Quotienten genau $(B: \nabla_{y} [\nu])$-mal 
als direkter Summand 
vor. Analoges gilt f"ur $\cal{F}_\Delta.$
\ELemma
\begin{proof}
Ist $z_{0}, z_{1}, z_{2}, \ldots $ eine Aufz"ahlung der Elemente
von $\CW$ derart, da"s die Folge $l(z_{0}), l(z_{1}), l(z_{2}),
\ldots$ der L"angen monoton w"achst, und setzen wir $A(j) = \{
z_{0}, z_{1}, \ldots , z_{j}\}$, so bilden die $\Gamma_{A(j)} B$
eine Verfeinerung unserer Filtrierung $\Gamma_{\leq i}B$ und
$\Gamma_{A(j)}B /\Gamma_{A(j-1)}B$ ist eine Summe von 
$(B:\nabla_{z_j}[\nu])$ Kopien von
gewissen $\nabla_{z_j}[\nu].$ Zwischen je zwei  
Aufz"ahlungen der Elemente
von $\CW$ derart, da"s in der Bruhat-Ordnung gr"o"sere
Elemente stets auch einen gr"o"seren Index haben, k"onnen wir jedoch
offensichtlich hin- und hergehen in endlich vielen Schritten, 
von denen jeder nur zwei benachbarte unvergleichbare Elemente
vertauscht. Da sich unvergleichbare Elemente nicht um
eine Spiegelung unterscheiden k"onnen, gibt es bei jedem Schritt
kein $\EExt^{1}$ zwischen den entsprechenden Subquotienten, 
folglich liefern die Filtrierungen vor und nach jedem dieser
Schritte bis auf Reihenfolge dieselben Subquotienten.
Der Subquotient
$\Gamma_{C(k)} B/\Gamma_{C(k-1)} B$ mit $y_{k} =z_{j}$ ist
also isomorph zu
$\Gamma_{A(j)}B/\Gamma_{A(j-1)} B,$
mithin ist $\Gamma_{\leq y} B/ \Gamma_{<y} B$ eine Summe von
verschobenen Kopien von $\nabla_{y}$, und zwar kommt $\nabla_{y}
[\nu]$ genau $(B: \nabla_{y} [\nu])$-mal vor.
\end{proof}

\begin{Proposition}\label{frei}
Gegeben $B\in\CB$ und $y\in \CW$ sind die eben definierten
$\Gamma^\leq_y  B,$ 
$\Gamma^\geq_y  B$ und
$\Gamma^y B$ ebenso wie  $\Gamma_y B$ graduiert freie $R$-Rechtsmoduln,
auf denen die Operation von $R\otimes R$  faktorisiert "uber $R_y.$
\end{Proposition}
\begin{proof}
F"ur die ersten beiden Moduln unserer Liste folgt das
leicht aus \ref{PP}.
F"ur den Letzten folgt dann aus der Einbettung
$\Gamma_y B\hra \Gamma^\leq_y  B$ nach \ref{BB},
da"s die Operation von $R\otimes R$  "uber $R_y$ faktorisiert, 
da"s also das Auswerten bei $1_y$ einen 
Isomorphismus $\op{Hom}(R_y,B)\sira\Gamma_y B$ liefert,
der mit 
$\ref{DHOM}$ dann die Freiheit von $\Gamma_y B$ zeigt.
Die Behauptung f"ur $\Gamma^y B$ schlie"slich folgt
f"ur $y=e$  aus $B\in\CF_\Delta.$
Im allgemeinen beachten wir, da"s aus der Definition
der Kategorie $\CB$ und \ref{HD} sogar folgt 
$B\otimes_R R_z\in \CF_\Delta$ f"ur alle $z\in\CW.$
Wenden wir diese Erkenntnis an auf $z=y^{-1},$  so
ergibt sich der allgemeine Fall.
\end{proof}

\begin{Notation}
Bezeichne $\CT \subset \CW$ die Menge aller Spiegelungen, und f"ur
$t \in \CT$ bezeichne $\al_{t} \in V^{\ast}$ eine Gleichung seiner
Spiegelebene, in Formeln $\ker \al_{t} = V^{t}$.
Die $\al_{t}$ sind eindeutig bis auf einen Skalar, wir w"ahlen sie
f"ur den Rest dieses Abschnitts beliebig aber fest.
F"ur $y \in \CW$ betrachten wir nun in $R$ das Element
$$p_{y} = \prod_{t \in \CT\!,\; yt < y} \al_{t}.$$  
\end{Notation}

\begin{Satz}\label{IP}
Sei $y \in \CW$  und $B \in \CB$.
So induzieren die offensichtlichen Morphismen
 $\Gamma_{y}B\ra \Gamma^{\leq}_{y} B$
und $\Gamma_y^{\geq}B\ra \Gamma^{y}B$ Isomorphismen
$$
\begin{array}{lll}
\Gamma_{y}B&\sira& \Gamma^\leq_y  Bp_{y},\\[2mm]
\Gamma_{y}^{\geq}B&\sira& \Gamma^y  Bp_{y}.
\end{array}$$
\end{Satz}

\BBemerkung\label{CX}
Der Beweis wird erst im Anschlu"s an \ref{LS} gegeben.
In \cite{HCH} wird f"ur den Fall 
$|\CW|<\infty$ gezeigt, da"s es Elemente 
$c_{y} \in R
\otimes R$ vom Grad $2l(y)$ gibt,
die verschwinden auf $\Gr(x)$ f"ur $x < y,$ jedoch nicht auf $\Gr(y)$
selbst. 
Dann stimmt notwendig $c_y$
auf $\Gr(y)$ bis auf
einen Skalar mit $1\otimes p_y$ "uberein und wir folgern ein 
bis auf den fraglichen Skalar kommutatives
Diagramm von Isomorphismen
$$\begin{array}{lcc}
 \Gamma^\leq_y  B&\stackrel{c_y}{\ra}&\Gamma_{y}B\\[2mm]
{\scriptscriptstyle \cdot p_y }\da\;\;&&\parallel\\[2mm]
\Gamma^\leq_y  Bp_{y}&\leftarrow& \Gamma_{y} B
\end{array}$$
das das Inverse des ersten im Satz behaupteten
Isomorphismus etwas expliziter angibt.
"Ahnlich erhalten wir auch ein bis auf den fraglichen Skalar kommutatives
Diagramm mit Isomorphismen im Quadrat ganz rechts
$$\begin{array}{ccccccc}
\Gamma_{\neq y}B&\hra&B&\sra&\Gamma^y B&\sira&\Gamma^y Bp_y\\
\da&&\da&&\da&&\da\\
\Gamma_{\not\leq y}B&\hra&
\Gamma_{\not< y}B&\sra&\coker &\cong&\Gamma^{\geq}_y B\end{array}$$
wo alle Vertikalen bis auf die Vertikale
ganz rechts durch Multiplikation mit $c_y$
definiert sind, die obere Horizontale ganz rechts die Multiplikation mit
$p_y$ ist, und wir f"ur die
letzte untere Horizontale ein Analogon von
\ref{PP} verwenden, um eine kurze exakte Sequenz
$\Gamma_{\not\leq y}B\hra\Gamma_{\not< y}B\sra\Gamma^{\geq}_y B$
herzuleiten. So kann man das  Inverse des zweiten im Satz behaupteten
Isomorphismus etwas expliziter verstehen.
\end{Bemerkung}
\begin{Bemerkung}
Ich w"u"ste gerne, ob auch f"ur jedes Element $y$ einer
unendlichen Coxetergruppe
diejenige Funktion  $c_{y}:\op{Gr}(\leq y)\ra \DC$ regul"ar ist,
die verschwindet auf $\Gr(x)$ f"ur $x < y$ und die auf $\Gr(y)$
mit $1\otimes p_y$ "ubereinstimmt. 
\end{Bemerkung}
\begin{Notation}
F"ur eine Spiegelung $t\in \CW$ bezeichne $R_{(t)}$ 
die Lokalisierung von $R$ an allen
homogenen Funktionen auf $V$, die auf der Spiegelebene
$V^{t}$ nicht identisch verschwinden.
\end{Notation}
\begin{Lemma}\label{LS}
F"ur jeden speziellen Bimodul
$B\in \CB$ ist seine Lokalisierung 
$B \otimes_{R} R_{(t)}$ in $R\op{-mod_\DZ-} R_{(t)}$
isomorph zu einem direkten Summanden in einer endlichen direkten Summe von 
Kopien von $R_{y,yt} 
\otimes_{R} R_{(t)}$ und $R_{y} \otimes_{R}
R_{(t)}$ f"ur $y \in \CW$ mit $y < yt$.
\end{Lemma}
\BBemerkung
Die etwas m"uhsame Formulierung erlaubt es, eine Diskussion von
Krull-Schmid in diesem Zusammenhang zu vermeiden.
\EBemerkung
\begin{proof}
Wir zeigen das durch Induktion in Anlehnung an die induktive
Definition der Objekte von $\CB.$
Es reicht damit zu zeigen, da"s f"ur jede einfache Spiegelung $s \in
\CS$ die beiden $R$-$R_{(t)}$-Bimoduln
$R\otimes_{R^{s}}R_{y}\otimes_{R} R_{(t)}$ und
$R\otimes_{R^{s}} R_{y,yt} \otimes_{R} R_{(t)}$
eine Zerlegung in Bimoduln derselben Gestalt haben.
Nun k"onnen aber $R_{y}\otimes_{R}R_{(t)}$ und $R_{x}\otimes_{R}
R_{(t)}$ nur dann erweitern in $R\op{-mod_\DZ-} R_{(t)},$ wenn gilt $y
=xt.$ In der Tat gibt es nach 
\ref{Ext} bereits in $R\op{-mod_\DZ-} R$ nur dann
Erweiterungen, wenn gilt $y=xr$ f"ur eine Spiegelung $r,$ und
diese Erweiterungen sterben unter der Multiplikation von rechts
mit jeder Gleichung der Spiegelebene von $r.$

Um den ersten Fall zu erledigen m"ussen wir also nur zeigen, da"s
aus $sy \neq yt$ folgt $syt > sy$. Das folgt jedoch
aus der Deodhar's Eigenschaft Z oder alternativ aus der
geometrischen Erkenntnis,
da"s benachbarte Alkoven eben nur von einer Wand
getrennt werden.
Im zweiten Fall ist nichts zu tun 
falls gilt $sy =yt$, in diesem speziellen 
Fall
gilt
ja sogar $R\otimes_{R^{s}} R_{y,yt} \cong R_{y,yt} \oplus R_{y,yt}[-2]$
nach \ref{TTT}.

Gilt schlie"slich $sy \neq yt$, so f"uhrt die kurze exakte Sequenz
von Bimoduln
$R[-2]\hra R\otimes_{R^s}R\sra R_s$ zu einer kurzen exakten Sequenz
$$R_{y,yt}[-2]\hra R\otimes_{R^{s}} R_{y,yt}\sra R_{sy,sy t}$$
und diese Sequenz mu"s spalten nach Anwenden von $\otimes_{R} R_{(t)},$ 
da alle Erweiterungen zwischen den Bimoduln 
$R_{y},$ $ R_{yt}$ als Untermoduln
und $R_{sy},$ $ R_{syt}$ als Quotienten
spalten nach Anwenden von $\otimes_{R} R_{(t)}.$
\end{proof}\noindent
\begin{proof}[Beweis von \ref{IP}]
Der Kokern der Einbettung 
$\Gamma_{y}B\hra \Gamma_{\leq y} B$ enth"alt keine Elemente mit
Tr"ager $\op{Gr}(y).$  
Folglich induziert  die
Restriktion eine Einbettung
$$\HHom (\Gamma^{\leq}_{y} B, R_{y})=
\HHom (\Gamma_{\leq y} B, R_{y}) \hookrightarrow \HHom
(\Gamma_{y} B, R_{y}).$$
Wir zeigen zun"achst, da"s diese Einbettung sogar im Teilraum
$\HHom (\Gamma_{y} B, R_{y}p_{y})$ landet.
Gehen wir n"amlich f"ur eine Spiegelung $t\in \CT$ mit $yt < y$ durch
$\otimes_{R}R_{(t)}$ zu den Lokalisierungen "uber
und zerlegen $B$ wie im Lemma \ref{LS}, so tragen nur die
Summanden der Gestalt $R_{y,yt}\otimes_R R_{(t)}$ zu 
den beiden zu vergleichenden Hom-R"aumen bei
und wir folgern, da"s
unsere Einbettung landet in $\HHom (\Gamma_{y}B,
R_{y} \al_{t} \otimes_{R} R_{(t)}),$ 
weil das eben gilt, wenn wir
$B$ durch  $R_{y,yt}$ ersetzen.
Das gilt nun f"ur alle $t$ mit $yt
< y,$ und da der Schnitt aller dieser
$R_{y} \al_{t} \otimes_{R} R_{(t)}$ mit $R_y$ gerade $R_{y}p_y$ ist,
ergibt sich die gew"unschte Faktorisierung
$$\HHom (\Gamma^{\leq}_{y} B, R_{y})=
\HHom (\Gamma_{\leq y} B, R_{y}) \hookrightarrow \HHom
(\Gamma_{y} B, R_{y} p_{y}).$$
Da $\Gamma^{\leq}_{y} B$ nach \ref{frei} ein freier $R_y$-Modul ist,
folgern wir weiter, da"s die Einbettung $\Gamma_{y}B\hra \Gamma^{\leq}_{y} B$
bereits in $\Gamma^{\leq}_{y} Bp_y$ landet.
Wir zeigen nun durch Dimensionsvergleich,
da"s die so erhaltene Einbettung  ein 
Isomorphismus
ist. Lemma \ref{PP} zeigt schon
$$\underline{\rk} \Gamma^\leq_y  B=\sum (B: \nabla_{y}
[\nu])v^{\nu+l(y)}.$$
Auf der anderen Seite faktorisiert nach \ref{BB} die Operation
von $R\otimes R$ auf $\Gamma_{y} B$ "uber $R_y$ und wir haben folglich
$\Gamma_{y} B = \HHom (R_{y},B) = \HHom (\Delta_{y} [l(y)], B).$
Letzteres ist jedoch nach 
\ref{DHOM} 
ein freier $R$-Rechtsmodul vom Rang $\underline{\rk} \Gamma_{y} B = \sum (B:
\nabla_{y} [\mu]) v^{\mu -l(y)}.$ Da aber nun
$p_{y}$ im Grad $2l(y)$ lebt, stimmt das "uberein mit dem
Rang von $\Gamma^\leq_y  Bp_y$ 
und wir erhalten den ersten Isomorphismus in unserem Satz.

Um  den zweiten Isomorphismus herzuleiten
beachten wir dual, da"s der Kokern der Einbettung $\Gamma_{\geq y}B\hra B$
keine Elemente mit Tr"ager $\op{Gr}(y)$ besitzt. Folglich induziert
die Restriktion eine Einbettung
$$\HHom(B,R_y)\hra \HHom(\Gamma_{\geq y}B,R_y)=
 \HHom(\Gamma^\geq_y B,R_y ).$$
Wir zeigen zun"achst wie eben mit Lemma \ref{LS} "uber die lokale Struktur der
Bimoduln aus $\CB,$
da"s diese Einbettung sogar im Teilraum 
 $\HHom(\Gamma_{\geq y}B,R_yp_y)$ landet.
Dann zeigt ein Dimensionsvergleich der homogenen Teile, da"s 
die so definierte Einbettung ein Isomorphismus 
$$\HHom(B,R_y)\sira \HHom(\Gamma^\geq_y B,R_y p_y)$$
ist.
Daraus erhalten wir 
je nach Geschmack durch Dimensionsvergleich
oder Dualisieren  mit \ref{frei} die zweite Behauptung des Satzes.
\end{proof}

\begin{Bemerkung}
Wir werden bald zeigen k"onnen, da"s $\CB$ stabil ist unter 
dem Bilden von direkten Summanden. Solange wir das jedoch noch nicht wissen,
definieren wir $\add\CB$ als die Kategorie aller graduierten
$R$-Bimoduln, die sich als direkte Summanden von speziellen Bimoduln
realisieren lassen. Nun verallgemeinern wir \ref{HOM} 
auf $\add\CB$ und zeigen  
\end{Bemerkung}

\begin{Lemma}\label{HOOM}
F"ur $M \in \CF_{\Delta}$, $ N \in \add\CB$ 
und desgleichen f"ur $M\in\add\CB, $ $N\in \CF_\nabla$ ist $\HHom
(M,N)$ graduiert frei als $R$-Rechtsmodul vom Rang 
$$\underline{\rk}\HHom
(M,N)=
\sum_{x,\nu,\mu} (M: \Delta_{x} [\nu])(N:\nabla_{x} [\mu])v^{\mu
-\nu}.$$
\end{Lemma}
\begin{proof}
Wir beginnen mit dem ersten Fall.
Zun"achst beachten wir, da"s 
Isomorphismus $\Gamma_y N\sira \Gamma_y^{\leq} Np_y$ aus
Satz \ref{IP} 
auch f"ur alle $N\in\add\CB$ gelten mu"s.
Daraus folgt  die Aussage des Lemmas f"ur
$M=\Delta_y$ und $N\in\add \CB,$ denn
wir haben $$
\begin{array}{cl}
\underline{\rk}\op{Hom}(\Delta_y, N)
&=\underline{\rk}\Gamma_y N[l(y)]\\[2mm]
&=\sum_\mu(\Gamma_y N[l(y)]:\nabla_y[\mu])v^{\mu+l(y)}\\[2mm]
&=\sum_\mu(\Gamma_y^{\leq}N[-l(y)]:\nabla_y[\mu])v^{\mu+l(y)}\\[2mm]
&=\sum_\mu(N:\nabla_y[\mu])v^{\mu}.
\end{array}
$$
Im allgemeinen argumentieren
wir dann mit Induktion "uber die L"ange einer $\Delta$-Fahne von $M.$
Sei in der Tat $x\in\CW$ ein Element
maximaler L"ange mit $(M : \Delta_{x} [\nu]) \neq 0$ f"ur ein $\nu
\in \DZ$.
Wegen der Maximalit"at von $x$ haben wir eine kurze exakte Sequenz
$$\Gamma_{x} M \hookrightarrow M \twoheadrightarrow \coker$$
mit $\Gamma_{x}M$ und $ \coker $ in $ \CF_{\Delta}$ und
$$(M:\Delta_{y}[\nu]) = (\Gamma_{x}M : \Delta_{y}[\nu]) +(\coker :
\Delta_{y}[\nu])$$ f"ur alle $y \in \CW$ und $\nu \in \DZ$.
Nat"urlich ist hier sogar auf der rechten Seite der erste Summand Null
im Fall $y \neq x$ und der zweite im Fall $y =x$.
In jedem Fall aber erzwingen unsere Dimensionsformeln \ref{HOM} f"ur alle $N 
\in
\CB$ eine kurze exakte Sequenz
$$\HHom (\coker, N) \hookrightarrow \HHom (M,N) \twoheadrightarrow
\HHom (\Gamma_{x}M, N).$$
Diese Sequenz mu"s dann  auch f"ur alle $N\in\add\CB$ exakt sein,
und so zeigt unsere Induktion "uber die L"ange einer
$\Delta$-Fahne von $M$ umgekehrt die behaupteten Dimensionsformeln 
f"ur alle $N\in\add\CB.$

Den zweiten Fall behandeln wir analog. 
Da"s \ref{IP} auch f"ur alle $B\in\add \CB$ gelten mu"s
erledigt die F"alle $N=\nabla_y.$ Im allgemeinen
argumentieren wir "uber die L"ange einer $\nabla$-Fahne
von $N,$ 
nehmen  ein Element $x$ maximaler L"ange mit
$(N: \nabla_x[\nu])\neq 0$ f"ur mindestens ein $\nu\in\DZ,$ 
haben wegen der Maximalit"at von $x$ eine kurze exakte Sequenz
$$\ker \hookrightarrow N \twoheadrightarrow \Gamma^{x}N$$
mit $\Gamma^{x}N $ und $\ker$ in $ \CF_{\nabla}$ und $$(N : \nabla_{y}[\nu]
) = (\ker : \nabla_{y}[\nu])+(\Gamma^{x}N : \nabla_{y}[\nu])$$
f"ur alle $y \in \CW$ und $\nu \in \DZ$, und wieder erzwingen
unsere Dimensionsformeln \ref{HOOM} 
f"ur alle $M \in \CB$ eine kurze exakte
Sequenz
$$\HHom (M, \ker)\hookrightarrow \HHom (M, N)\twoheadrightarrow
\HHom (M, \Gamma^{x}N),$$
die dann auch f"ur alle $M\in \add\CB$ exakt sein mu"s und den
gew"unschten Induktionsbeweis erlaubt.
\end{proof}
\begin{Satz}\label{KL}\label{DSS}
\begin{enumerate}
\item
F"ur alle $x \in \CW$ gibt es bis auf Isomorphismus genau einen
unzerlegbaren speziellen Bimodul $B_{x} \in \CB$ mit Tr"ager in
$\Gr ( \leq x)$ und $(B_{x} : \Delta_{x}[\nu]) =1$ f"ur $\nu=0$ und
Null f"ur $\nu\neq 0.$
\item
Die Abbildung $(x,\nu)\mapsto B_x[\nu]$ definiert eine Bijektion
$$\CW\times\DZ\sira\left\{\begin{array}{c}\text{Unzerlegbare Objekte in 
$\CB,$}\\
\text{bis auf Isomorphismus}\end{array}\right\}$$
\item
Die Bimoduln $B_{x}$ sind selbstdual, in Formeln $DB_{x} \cong B_{x}$.
\item
Die Kategorie $\CB$ ist stabil unter dem Bilden von direkten Summanden, 
in Formeln $\CB=\add\CB.$
\end{enumerate}
\end{Satz}
\begin{proof}
Wir zeigen zun"achst die drei ersten Punkte f"ur $\add\CB$ statt
$\CB.$
Sei $x = st \ldots r$ eine reduzierte Darstellung.
In $R\otimes_{R^{s}} R \otimes_{R^{t}} \ldots R
\otimes_{R^{r}} R [l(x)],$  
das ja nach \ref{ZD} bis auf die Unzerlegbarkeit bereits
alle von $B_x$ geforderten Eigenschaften hat,
betrachten wir 
den unzerlegbaren Summanden $B_{x}$ mit
$(B_{x} : \Delta_{x}) =1$ und haben damit schon ein m"ogliches
$B_{x}\in\add\CB$ gefunden, das sogar offensichtlich auch noch selbstdual
ist als der einzige Summand eines 
selbstdualen Bimoduls, dessen Tr"ager $\op{Gr}(x)$ umfa"st.

Es gilt nun nur noch zu zeigen, da"s es au"ser den eben konstruierten
$B_{x}$ und ihren Verschobenen $B_{x}[\nu]$ keine unzerlegbaren
Objekte in $\add\CB$ gibt.
Aber sei $M \in \add\CB$ und sei $x \in \CW$ ein Element
maximaler L"ange mit $(M : \Delta_{x} [\nu]) \neq 0$ f"ur ein $\nu
\in \DZ$.
Wir zeigen, da"s $B_x[\nu]$ ein direkter Summand ist von $M.$
Wegen der Maximalit"at von $x$ liefern zun"achst nach
\ref{HOOM} oder genauer seinem Beweis die offensichtlichen 
Abbildungen Surjektionen
$$\begin{array}{l}
\HHom (M,B_x) \twoheadrightarrow
\HHom (\Gamma_{x}M, B_x)=\HHom (\Gamma_{x}M, \Gamma_{x}B_x),\\[2mm]
\HHom (B_x, M)\twoheadrightarrow
\HHom (B_x, \Gamma^{x}M)=\HHom (\Gamma^{x}B_x, \Gamma^{x}M).
\end{array}$$
Schlie"slich beachten wir noch, da"s aufgrund der Maximalit"at von
$x$ mit \ref{PP} sogar gilt $\Gamma_{x}^{\leq} M \sira \Gamma^{x}M$, so da"s
also aufgrund von \ref{IP} die Einbettung $\Gamma_{x} M
\hookrightarrow \Gamma^{x}M$  Isomorphismen $$i : \Gamma_{x} M
\sira \Gamma^{x} M p_{x} \quad\text{und}\quad i : \Gamma_{x} B_x
\sira \Gamma^{x} B_x p_{x}$$ induziert.
Gegeben $m \in
\Gamma^{x}M$ homogen finden wir nun $f \in \HHom (\Gamma^{x}B_{x},
\Gamma^{x}M)$ homogen mit $f(b) = m$ f"ur $b \in \Gamma^{x}B_{x}$
einen homogenen Erzeuger.
W"ahlen wir $m$ sogar so, da"s es sich zu einer homogenen Basis
des graduiert freien $R$-Moduls $\Gamma^{x}M$ erg"anzen l"a"st, so
finden wir auch $g \in \HHom (\Gamma_{x}M, \Gamma_{x}B_{x})$
mit $g(i^{-1}(mp_{x})) = i^{-1}(bp_{x})$.
F"ur homogene Lifts $\tilde{f} : B_{x} \ra
M$ und $\tilde{g} : M \ra B_{x}$
von $f$ und $g$
folgt, da"s $\tilde{g} \circ
\tilde{f} : B_{x} \ra B_{x}$ auf $\Gamma^{x}B_{x}$ die Identit"at
induziert.
Also kann $\tilde{g} \circ \tilde{f}$ nicht nilpotent sein, also
ist es ein Isomorphismus, da ja $B_{x}$ unzerlegbar ist.
Schlie"slich folgt 
durch Induktion "uber die L"ange von $x,$ 
da"s alle $B_x$ schon zu $\CB$ geh"oren:
In der Tat zerf"allt unser langes Tensorprodukt oben,
das per definitionem zu $\CB$ geh"ort,
in einen Summanden $B_x$ nebst anderen Summanden,
die wegen ihres Tr"agers die Gestalt $B_y[\nu]$ mit $l(y)<l(x)$ 
haben m"ussen.
Das zeigt die letzte Behauptung.
\end{proof}
\BBemerkung\label{BBn}
Wenn wir ein
$B \in\CB$ finden mit $\langle B \rangle=\CE(\SDu_x),$
so folgt schon $B\cong B_x$ und damit $\langle B_x
\rangle=\CE(\SDu_x).$ In der Tat ist ein solches $B$ notwendig
unzerlegbar, denn berechnen wir die Dimensionen der homogenen
Komponenten seines Endomorphismenrings, so finden wir nichts
in negativen Graden und nur den Grundk"orper im Grad Null.
\EBemerkung

\begin{Bemerkung}
F"ur alle $x \in \CW$ ist $h_{\Delta}(B_{x})= h_{\nabla}(B_{x})$
jedenfalls ein selbstduales Element von $\CH$ der Gestalt $\SDu_{x}
+ \sum_{y<x} h_{y}\SDu_{y}$ f"ur geeignete selbstduale 
$h_{y} \in \DZ[v, v^{-1}].$
Sei in der Tat $x = st \ldots r$ eine reduzierte Zerlegung.
So ist 
$$R\otimes_{R^{s}}R\otimes_{R^{t}} \ldots R\otimes_{R^{r}} R[l(x)]$$ 
selbstdual und zerf"allt folglich in $B_{x}$ plus
Summanden der Gestalt $B_{y}$ und $B_{y} [\nu] \oplus B_{y}[-\nu]$
f"ur $y<x$.
Wir folgern
$h_{\Delta} (B_{x}) = h_{\nabla}(B_{x})$ selbstdual mit Induktion
"uber die Bruhat-Ordnung, und dann gilt offensichtlich
$$h_{\Delta} (B_{x}) = h_{\nabla}(B_{x})=\SDu_{x}+
\sum_{y<x}h_{y} \SDu_{y}$$ f"ur geeignete 
selbstduale $h_{y}\in \DZ[v, 
v^{-1}]$.
Ich kann allerdings noch nicht einmal zeigen, da"s diese $h_{y}$ nichtnegative
Koeffizienten haben.  
\end{Bemerkung}

\section{Diskussion der Hauptvermutung}
Ein nat"urlicher Ansatz zum Beweis der Vermutung \ref{HV} scheint
mir voll\-st"an\-dige Induktion "uber die L"ange von $x.$
Ich formuliere im folgenden verschiedene Varianten dieses Ansatzes
in einem Omnibus-Lemma.
Bezeichne $\underline{\HHom}(B,B')$ die $(R \otimes
R)$-Homomorphismen vom $\DZ$-Grad Null
zwischen zwei $\DZ$-graduierten $R$-Bimoduln $B,B'.$
\BLemma 
Seien $x\in\CW$ und $s\in\CS$ gegeben mit $sx> x.$
Wir nehmen an, es gelte
$\langle B_{y}\rangle =
\CE(\SDu_{y})$ f"ur alle $y\in\CW$ mit $l(y)\leq l(x).$ 
Unter dieser Annahme sind gleichbedeutend:
\begin{enumerate}
\item
Es gilt $\langle B_{sx}\rangle =
\CE(\SDu_{sx});$
\item
F"ur alle $y$ mit
$l(y)\leq l(x)$ definiert die
Komposition von Morphismen 
$\underline{\HHom} (B_{y}, \theta_s B_{x})\times \underline{\HHom}
(\theta_s B_{x}, B_{y}) \ra \underline{\EEnd} B_{y} = {{k}}$
eine nicht ausgeartete Paarung von ${{k}}$-Vektorr"aumen.
\item
F"ur alle $y$ mit
$l(y)\leq l(x)$ definiert  
$\Gamma_y\theta_s B_x\hra \Gamma^y\theta_s B_x$ eine Inklusion 
$(\Gamma_y\theta_s B_x)_{l(y)}\hra \Gamma^y\theta_s B_x p_y \otimes_R{{k}}.$
\item
F"ur alle $y$ mit
$l(y)\leq l(x)$ definiert die offensichtliche Abbildung 
eine Inklusion 
$(\Gamma_y\theta_s B_x)_{l(y)}\hra \Gamma^\geq_y\theta_s B_x \otimes_R{{k}}.$
\item
F"ur alle $y$ mit
$l(y)\leq l(x)$ definiert die offensichtliche Abbildung 
eine Inklusion 
$(\Gamma^\leq_y\theta_s B_x)_{-l(y)}\hra \Gamma^y\theta_s B_x
\otimes_R{{k}}.$
\end{enumerate}
\ELemma
\begin{proof}
$1\RA 2:$ Dem Leser "uberlassen.
$2\RA 1:$ 
In der Notation von \cite{So-K}
k"onnen wir schreiben $$\SDu_{s} \SDu_{x} =
\SDu_{sx} + \sum_{{y < sx}} m_{y}
\SDu_{y}$$
mit gewissen $m_{y} \in \DZ$ und es gilt
$\CE(\SDu_{s} \SDu_{x})=\langle\theta_s B_{x}\rangle.$
Nach \ref{HOM} und Induktionsannahme haben wir also
$$
\dim \underline{\HHom} (B_{y}, \theta_s B_{x}) =m_y
 = \dim \underline{\HHom} (\theta_s B_{x}, B_{y}).
 $$
Ist unsere Paarung nicht ausgeartet, so k"onnen wir 
f"ur jedes $y$ jeweils $m_{y}$
Kopien von $B_{y}$ von $\theta_s B_{x}$ abspalten und so 
ein $B\in\CB$ konstruieren mit 
$\langle B\rangle =
\CE(\SDu_{sx}).$
Nach Bemerkung \ref{BBn} gilt dann aber notwendig 
$B\cong B_{sx}.$ 

$3\IFF 2:$ 
Mit Dimensionsvergleich erkennen wir, da"s 
auf den Homomorphismen vom Grad Null die  
offensichtlichen Abbildungen Isomorphismen 
$$\begin{array}{ccc}\underline{\HHom} (B_{y}, \theta_s B_{x}) 
& \sira & \underline{\HHom} 
(\Delta_y, \theta_s B_{x})\\
\underline{\HHom}(\theta_s B_{x}, B_{y}) &\sira & \underline{\HHom} 
(\theta_s B_{x}, \nabla_{y})
\end{array}$$
liefern.
Nach \ref{IP} wissen wir, da"s die Komposition der Morphismen
vom Grad Null
$\Delta_{y}\ra B_y \ra \nabla_{y}$ bis auf einen Skalar
die Rechtsmultiplikation mit $p_y$ ist, und da"s ganz allgemein 
f"ur alle $B\in\CB$ jede Komposition 
$\Delta_{y}\ra B \ra \nabla_{y}$ in $\nabla_{y}p_y$
landet.
Damit k"onnen wir unsere Paarung aus Teil 2 umschreiben zu einer Paarung
$$\underline{\HHom} (\Delta_{y}, \theta_s B_{x})\times \underline{\HHom}
(\theta_s B_{x}, \nabla_{y}) \ra p_y{{k}}\subset
\underline{\HHom} (\Delta_{y},\nabla_y).$$
Nun ist aber so eine Paarung zwischen zwei R"aumen nichts anderes als eine
Abbildung vom einen in den Dualraum des anderen,
wir identifizieren m"uhelos
$$\begin{array}{ccl}
\underline{\HHom} 
(\Delta_y, \theta_s B_{x})&=&(\Gamma_y\theta_s B_x)_{l(y)}\\[2mm]
\underline{\HHom}
(\theta_s B_{x}, \nabla_{y}) &=&
\underline{\HHom}
(\Gamma^y\theta_s B_x, \nabla_{y})\\
&=&
\underline{\HHom}
(\Gamma^y\theta_s B_x, {{k}}[l(y)])\\
&=&
\underline{\HHom}
((\Gamma^y\theta_s B_x)\otimes_R{{k}}, {{k}}[l(y)])\\[2mm]
\underline{\HHom}
(\theta_s B_{x}, \nabla_{y})^\ast &=&
((\Gamma^y\theta_s B_x)\otimes_R{{k}})_{l(y)}
\end{array}$$
und erkennen so $3\IFF 2.$
Die Varianten 4 und 5 sind Umformulierungen von 3 unter
Zuhilfenahme von \ref{IP}.
\end{proof}
\begin{Bemerkung}
Es scheint mir nun naheliegend, die Einbettungen
$$\Gamma_{y}B_{x} \hookrightarrow \Gamma^{y}B_{x}$$
n"aher zu untersuchen.
Definieren wir eine Surjektion $R \twoheadrightarrow {{k}} [T]$
als die Restriktion auf die Gerade ${{k}}{\rho}$, so vermute ich,
da"s der Kokern von $(\Gamma_{y}B_{x}) \otimes_{R} {{k}} [T]
\hookrightarrow (\Gamma^{y}B_{x} )\otimes_{R} {{k}} [T]$ eine
direkte Summe sein sollte von Kopien von $({{k}} [T]/(T^{i+1}))
[i]$, wo wir wieder $\deg T=2$ genommen haben.
Das sollte sogar allgemeiner gelten f"ur alle Geraden, die nicht ganz in einer
Spiegelebene enthalten sind.  
\end{Bemerkung}


\providecommand{\bysame}{\leavevmode\hbox to3em{\hrulefill}\thinspace}
\providecommand{\MR}{\relax\ifhmode\unskip\space\fi MR }
\providecommand{\MRhref}[2]{%
  \href{http://www.ams.org/mathscinet-getitem?mr=#1}{#2}
}
\providecommand{\href}[2]{#2}

\end{document}